% AMSTeX
\documentstyle{amsppt}
\magnification=1200
\hsize=6.5truein
\vsize=8.9truein
\topmatter
\title On the Hodge Filtration of Hodge Modules
\endtitle
\author Morihiko Saito\endauthor
\affil RIMS Kyoto Univ., Kyoto 606-8502 Japan \endaffil
\endtopmatter
\tolerance=1000
\baselineskip=12pt
\def\ssbull{\raise.2ex\hbox{${\scriptscriptstyle\bullet}$}}
\def\cB{{\Cal B}}
\def\cD{{\Cal D}}
\def\cE{{\Cal E}}
\def\cH{{\Cal H}}
\def\cM{{\Cal M}}
\def\cO{{\Cal O}}
\def\bC{{\Bbb C}}
\def\bD{{\Bbb D}}
\def\bN{{\Bbb N}}
\def\bQ{{\Bbb Q}}
\def\bR{{\Bbb R}}
\def\bZ{{\Bbb Z}}
\def\Gr{\text{{\rm Gr}}}
\def\Im{\hbox{{\rm Im}}}
\def\Coker{\hbox{{\rm Coker}}}
\def\Ker{\hbox{{\rm Ker}}}
\def\Hom{\hbox{{\rm Hom}}}
\def\DR{\hbox{{\rm DR}}}
\def\IC{\hbox{{\rm IC}}}
\def\Var{\hbox{{\rm Var}}}
\def\Sol{\hbox{{\rm Sol}}}
\def\Sp{\hbox{{\rm Sp}}}

\def\supp{\hbox{{\rm supp}}\,}
\def\Sing{\hbox{{\rm Sing}}\,}
\def\codim{\text{\rm codim}}
\def\reg{\text{\rm reg}}
\def\can{\text{\rm can}}
\def\simto{\buildrel\sim\over\to}
\def\SameAuthor{\vrule height3pt depth-2.5pt width1cm}

\document
\noindent
Let
$ X $ be a complex manifold, and
$ Z $ an irreducible closed analytic subset.
We have the polarizable Hodge Module
$ \IC_{Z}\bQ^{H} $ whose underlying perverse sheaf is the intersection
complex
$ \IC_{Z}\bQ $.
See [16].
Let
$ (M,F) $ be its underlying filtered
$ \cD_{X} $-Module.
Then
$ M $ is the unique regular holonomic
$ \cD_{X} $-Module which corresponds to
$ \IC_{Z}\bC $ by
the Riemann-Hilbert correspondence [9] [13], and it is relatively easy
to determine $ M $ in some cases (for example, if
$ Z $ is a hypersurface with isolated singularities [28]).
However, the Hodge filtration
$ F $ on
$ M $ is a more delicate object, and it is not easy to describe
$ F $ explicitly, because we have to calculate the filtration
$ V $ to some extent.

More generally, assume
$ (M,F) $ underlies a polarizable Hodge Module
$ \cM $ with strict support
$ Z $ (see [16]).
Let
$ q = \min\{p \in \bZ : F_{p}M \ne 0\} $.
The {\it generating level} of
$ (M,F) $ is defined to be the minimal length of a part of the filtration
$ F $ which generates
$ (M,F) $ over
$ (\cD_{X},F) $ (i.e.,
$ (M,F) $ has generating level
$ \le r $ if
$ F_{q+r+i}M = F_{i}\cD_{X}F_{q+r}M $ for
$ i \ge 0 $,
and generating level
$ r $ if it has generating level
$ \le r $ but not
$ \le r-1) $.
Here the filtration
$ F $ on
$ \cD_{X} $ is by the order of differential operator.
Note that
$ (M,F) $ has generating level
$ \le r $ for
$ r \gg 0 $ at least locally, because
$ F $ is a good filtration (i.e.,
$ \Gr^{F}M $ is coherent over
$ \Gr^{F}\cD_{X}) $.
The notion of generating level measures in some sense the complexity of
the filtration
$ F $.
(See also Remark (ii) after (1.3).)
If
$ Z $ is a point so that
$ \cM $ corresponds to a Hodge structure, then the generating level of
$ (M,F) $ coincides with the level of the Hodge structure.

\medskip\noindent
{\bf 0.1.~Theorem.} {\it If the restriction of
$ \cM $ to a Zariski-open subset of
$ Z $ is a variation of Hodge structure with stalkwise level
$ r $,
then
$ (M,F) $ has generating level
$ \le r + m - 1 $,
where
$ m = \dim Z $.
In particular, if
$ (M,F) $ underlies
$ \IC_{Z}\bQ^{H} $,
then it has generating level
$ \le m - 1 $.
}

\medskip
Actually we prove a more general statement by induction on
$ m $ using resolution of singularity and the calculation of nearby cycles
in the normal crossing case [17].
See (1.3) below.
The above estimate of the generating level is optimal, because we can
determine the generating level in a special case (see (0.8) below for a
more general statement):

\medskip\noindent
{\bf 0.2.~Theorem.} {\it Assume
$ (M,F) $ underlies
$ \IC_{Z}\bQ^{H} $ and
$ Z $ is a cone of a smooth hypersurface of degree
$ d $.
Let
$ k_{1} \in \bZ $ such that
$ k_{1} < \dim X - (\dim X/d) - 1 \le k_{1} + 1 $.
Then
$ (M,F) $ has generating level
$ k_{1} $.
In particular, it has generating level
$ m - 1 $ if
$ d > \dim X = m + 1 $.
}

\medskip
This gives a counter example to Brylinski's conjecture [3], and, more
generally, to the following statement (which seems to have been
conjectured by some people):

Assume
$ (M,F) $ underlies a polarizable Hodge Module
$ \cM $ with strict support
$ Z $,
and let
$ Z_{0} $ be the smallest closed subspace of
$ Z $ such that the restriction of
$ \cM $ to its complement is a variation of Hodge structure (see [16]).
Then there should exist a reduced lattice
$ L $ of
$ M $ such that
$$
F_{p}M = \sum_{i} F_{i}\cD_{X}(L \cap \tilde{j}_{*}(F_{p-i}
M|_{X\setminus Z_{0}})),
\leqno(0.3)
$$
where
$ \tilde{j} : X \setminus Z_{0} \to X $ is a natural inclusion.
(Here a lattice means an
$ \cO_{X} $-sub-Module of
$ M $ which generates
$ M $ over
$ \cD_{X} $,
and a lattice is called reduced if it is annihilated by the reduced
ideal of
$ \supp M $.)

In fact, if the variation of Hodge structure is of type
$ (0,0), $ then (0.3) for
$ p = q $ implies
$ L|_{X\setminus Z_{0}} = F_{q}M|_{X\setminus Z_{0}} $ where
$ q = \codim_{X} Z $ (using (1.5.3) below).
So (0.3) becomes
$ F_{p}M = F_{p-q}\cD_{X} L $,
i.e.,
$ F $ has generating level
$ 0 $.
This contradicts (0.2).

For the proof of (0.2), assume first
$ Z $ is a hypersurface of
$ X $,
and let
$ f $ be a local (reduced) equation of
$ Z $.
We will denote
$ M $ by
$ M_{f} $ in this case.
Let
$ M_{f}' = \cO_{X}[f^{ -1}] $,
and
$ M_{f}'' = M_{f}'/\cO_{X} $.
Then
$ M_{f} $ is (up to a Tate twist) isomorphic to the minimal
$ \cD_{X} $-sub-Module of
$ M_{f}'' $ whose restriction to
$ X \,\setminus \, \Sing Z $ coincides with the restriction of
$ M_{f}'' $.
See for example [28].
(In fact, the assertion holds for the corresponding perverse sheaves [1].)
Furthermore,
$ M_{f}', M_{f}'' $ have the Hodge filtration
$ F $,
because they naturally underlie mixed Hodge Modules.
See (2.2).
Let
$ \tilde{V} $ be the filtration on
$ \cO_{X} $ induced by the filtration
$ V $ on the algebraic microlocalization
$ \cO_{X}[\partial_{t},\partial_{t}^{-1}] $ in
[19] using the natural inclusion
$ \cO_{X} \to \cO_{X}[\partial_{t},\partial_{t}^{-1}] $.
Here
$ V $ is indexed decreasingly by
$ \bQ $ so that the action of
$ \partial_{t}t - \alpha $ on
$ \Gr_{V}^{\alpha}(\cO_{X}[\partial_{t},\partial_{t}^{-1}]) $ is
nilpotent.
See (2.1).
Using the relation of the Hodge filtration with the filtration
$ V $ in (2.3.3), we can show (see (2.4)):

\medskip\noindent
{\bf 0.4.~Theorem.} {\it
$ F_{0}M_{f}' = f^{ -1}\tilde{V}^{1}\cO_{X} $,
$ F_{p}M_{f}' \supset \sum_{k\ge 0} F_{p-k}\cD_{X}
(f^{-k-1}\tilde{V}^{k+1}\cO_{X}) $.
}

\medskip
We have a similar formula for
$ M_{f}'' $,
because
$ (M_{f}',F) \to (M_{f}'',F) $ is strictly surjective (i.e.,
$ F $ on
$ M_{f}'' $ is the quotient filtration).
For example, we have
$ F_{0}M_{f}'' = f^{ -1}\tilde{V}^{1}\cO_{X}/\cO_{X} $.
Using the theory of duality for filtered differential complexes [16], we
can show that
$ \omega_{X} \otimes_{\cO_{X}} F_{0}M_{f}'' $ is the dual of
$ \Gr_{F}^{0}{\underline{\Omega}}_{Z} $,
where
$ (\underline{\Omega}_{Z},F) $ is the filtered complex of
du Bois [5].
We say that
$ Z $ has du Bois singularity if the natural morphism
$ \cO_{Z} \to \Gr_{F}^{0}{\underline{\Omega}}_{Z} $ is
a quasi-isomorphism [25].
Let
$ b_{f}(s) $ be
the
$ b $-function of
$ f $,
and
$ \alpha_{f} $ the minimal root of
$ b_{f}(-s)/(s - 1) $.
Then

\medskip\noindent
{\bf 0.5.~Theorem.} {\it
$ Z $ has
du Bois singularity if and only if
$ \alpha_{f} \ge 1 $ {\rm (}i.e.,
the maximal root of
$ b_{f}(s) $ is
$ - 1) $.
}

\medskip\noindent
In fact, the both conditions are equivalent to
$ \tilde{V}^{1}\cO_{X} = \cO_{X}, $ see (2.4).
For the isolated singularity case, see also [25].

Note that (0.4) is a refinement of [18, (0.11)] on the relation between
the Hodge and pole order filtrations, because the last term of (0.4)
coincides with
$ f^{ -p-1}\cO_{X} $ for
$ p \le \alpha_{f} - 1 $ where
$ \alpha_{f} $ is as in (0.5).
From (0.4) we can deduce (see (2.5)):

\medskip\noindent
{\bf 0.6.~Theorem.} {\it We have
$ F_{1}M_{f} = (f^{ -1}\tilde{V}^{>1}\cO_{X})/\cO_{X},
(f^{-k- 1}\tilde{V}^{>k+1}\cO_{X})/\cO_{X} \subset M_{f} $
and
$$
F_{p+1}M_{f} \supset \sum_{k\ge 0} F_{p-k}\cD_{X}(f^{ -k-1}
\tilde{V}^{>k+1}\cO_{X})/\cO_{X}.
$$
}

\medskip
Note that the first equality means
$ \omega_{X} \otimes_{\cO_{X}} (f^{ -1}\tilde{V}^{>1}
\cO_{X})/\cO_{X} = \pi_{*}\omega_{Z'} $ where
$ \pi : Z' \to Z $ is a resolution of singularity (see [18]).
So we have
$ \omega_{X} \otimes_{\cO_{X}} (\cO_{X}/\tilde{V}^{>1}
\cO_{X}) = \omega_{Z} /\pi_{*}\omega_{Z'} $, and
$ Z $ has rational singularity if and only if
$ \alpha_{f} > 1 $ (see [loc. cit.] and compare (0.5)).
It is not clear whether the equality should hold in the last formula of
(0.4), (0.6).
We can prove the equality in a special case as in (0.7), (0.8) below.
If
$ Z $ is smooth, we have
$ M_{f} = M_{f}'' $,
and
$ F_{p}M_{f}' = \cO_{X}f^{ -p-1}, F_{p}M_{f}'' = F_{p+1}M_{f} =
\cO_{X}f^{ -p-1}/\cO_{X} $ for
$ p \ge 0 $.
So we may restrict to a neighborhood of
$ \Sing Z $.

Assume further that
$ \Sing Z $ is isolated,
$ n := \dim X > 1 $ and
$ f $ is a quasihomogeneous polynomial of weights
$ (w_{1}, \dots, w_{n}) $ for a local coordinate system
$ (x_{1}, \dots, x_{n}) $ on a neighborhood of
$ 0 \in \Sing Z $.
Then we have the order function
$ v_{w} $ on
$ M_{f,0}' $ such that
$ v_{w}(f^{ k}m) = v_{w}(m) + k $ for
$ m \in M_{f,0}', k \in \bZ $,
and
$ v_{w}(1) = \sum_{i} w_{i} = \alpha_{f} $.
See (4.1).
Let
$ G $ be the filtration on
$ M_{f}', M_{f}'' $ by pole order (i.e.,
$ G_{k}M_{f}' = \cO_{X}f^{ -k-1} $ if
$ k \ge 0 $ and
$ 0 $ otherwise).
For
$ \beta \in \bQ $,
let
$ G_{k}^{\ge \beta}M_{f,0}', G_{k}^{>\beta}M_{f,0}' $ be the
subgroups of
$ G_{k}M_{f,0}' $ defined by the condition
$ v_{w}(m) \ge \beta $ and
$ v_{w}(m) > \beta $ respectively, where
$ M_{f,0}' $ is the stalk of
$ M_{f}' $ at
$ 0 $.
Let
$ G_{k}^{\ge \beta}M_{f,0}'', G_{k}^{>\beta}M_{f,0}'' $ be the
image of
$ G_{k}^{\ge \beta}M_{f,0}', G_{k}^{>\beta}M_{f,0}' $ in
$ M_{f,0}'' $.
Then we have (see (4.2)):

\medskip\noindent
{\bf 0.7.~Theorem.} {\it With the above assumptions and notation, let
$ k_{0} = [n - \alpha_{f}] - 1 $.
Then
$$
F_{p}M_{f,0}' = \sum_{k\ge 0} F_{p-k}\cD_{X,0}G_{k}^{\ge
0}M_{f,0}' = \sum_{0\le k\le k_{0}} F_{p-k}\cD_{X,0}G_{k}^{\ge
0}M_{f,0}',
$$
and
$ (M_{f}',F) $ has generating level
$ k_{0} $ {\rm (}similarly for
$ M_{f}'') $.
}

\medskip
In particular, the equality holds in the last formula of (0.4), because
the middle term of the formula of (0.7) is contained in the last term of
(0.4) (see (4.1.2) below).
For the proof of (0.7), we have to calculate the filtration
$ V $ along
$ f $ which determines the Hodge filtration
$ F $ on
$ M_{f}' $.
Using the algebraic microlocalization on which the filtration
$ V $ is much easier to calculate (see (4.1.2)), the assertion is reduced
to a certain calculation about the regular sequence
$ (f_{1}, \dots, f_{n}) $ where
$ f_{i} = \partial f/\partial x_{i} $.
The assertion for
$ M_{f}'' $ follows from that for
$ M_{f}' $.
From (0.7) we can deduce (see (4.3)):

\medskip\noindent
{\bf 0.8.~Theorem.} {\it With the assumptions and notation of
{\rm (0.7),} we have
$ G_{k}^{>0}M_{f,0}'' \subset M_{f,0} $,
and
$$
F_{p+1}M_{f,0} = \sum_{k\ge 0} F_{p-k}\cD_{X,0}
G_{k}^{>0}M_{f,0}'' = \sum_{0\le k\le k_{1}} F_{p-k}
\cD_{X,0}G_{k}^{>0}M_{f,0}'',
$$
where
$ k_{1} = \max\{k \in \bZ : k < n - \alpha_{f} - 1\} $.
Furthermore
$ (M_{f},F) $ has generating level
$ k_{1} $.
}

\medskip
In particular, the equality holds in the last formula of (0.6) using
(4.1.2).
If
$ f $ is homogeneous of degree
$ d $ (i.e.,
$ w_{i} = 1/d) $,
then the first equality of (0.8) is related to a modified version
(see [16]) of Brylinski's conjecture.
See Remark after (4.3).
(However I am not sure whether the modified version should be true in
general.)
Note that (0.2) is a special case of of (0.8), because
$ \alpha_{f} = n/d $.

For the passage from (0.4), (0.7) to (0.6), (0.8), we have to work on
the difference between
$ \tilde{V}^{>k+1}, G_{k}^{>0} $ and
$ \tilde{V}^{k+1}, G_{k}^{\ge 0} $.
For (0.6) we use the theory of microlocal filtration
$ V $ in [18] so that the assertion is reduced to the fact that
$ M_{f} $ has no nontrivial quotient whose support has dimension
$ < m $.
For (0.8) we use Brieskorn's module
$ H_{f}'' $ in [2] and the algebraic local cohomology
$ H_{[0]}^{n}\cO_{X} $ to get an algebraic version of [28].
Then the assertion is reduced to the injectivity of the action of
$ f $ on Brieskorn's module [22].

Finally, the above two theorems can be extended to the following (see
(5.4)):

\medskip\noindent
{\bf 0.9.~Theorem.} {\it Except for the assertions on generating levels,
the assertions of {\rm (0.7)} and {\rm (0.8)} are true also in the
semiquasihomogeneous isolated singularity case.
}

\medskip
In this case,
$ (M_{f}',F) $ and
$ (M_{f},F) $ have generating level
$ \le k_{0} $ and
$ \le k_{1} $ respectively, because of the last isomorphisms of (0.7),
(0.8).
However, they may have strictly smaller generating levels. See Remark (i)
after (5.4).
In particular, we get

\medskip\noindent
{\bf 0.10.~Remark.} The generating levels of
$ (M_{f},F), (M_{f}',F) $,
and
$ (M_{f}'',F) $ are not necessarily constant under a
$ \mu $-constant deformation of
$ f $.

In \S 1 we introduce the notion of generating level, and prove (0.1).
In \S 2 we review the theory of (microlocal) filtration
$ V $ to show (0.4), (0.6).
In \S 3 we recall some facts from the theory of Brieskorn's module and
$ b $-function in the isolated hypersurface singularity case.
In \S 4 we restrict further to the quasihomogeneous isolated singularity
case, and prove (0.7), (0.8).
In \S 5 we study the generalization of (0.7), (0.8) to the
semiquasihomogeneous isolated singularity case.

\bigskip\bigskip\centerline{{\bf 1. Generating Level}}

\bigskip\noindent
{\bf 1.1.}
Let
$ X $ be a complex manifold,
$ (M,F) $ a filtered
$ \cD_{X} $-Module,
and
$ r $ an integer.
We say that
$ (M,F) $ is
$ r $-{\it generated} if
$$
F_{i+r}M = F_{i}\cD_{X}F_{r}M\quad
\text{for}\,\,i \ge 0,
\leqno(1.1.1)
$$
where
$ F $ on
$ \cD_{X} $ is the filtration by the order of differential operator.
We say that
$ (M,F) $ is exactly
$ r $-generated if it is
$ r $-generated but not
$ (r-1) $-generated.

Assume
$ q := \min\{p \in \bZ : F_{p}M \ne 0\} $ is finite.
We say that
$ (M,F) $ has {\it generating level}
$ \le r $ if it is
$ (r+q) $-generated, and generating level
$ r $ if it has generating level
$ \le r $ but not
$ \le r - 1 $.

\medskip\noindent
{\bf Remarks.} (i) We say that a filtered
$ \cD_{X} $-Module
$ (M,F) $ is coherent and the filtration
$ F $ on
$ M $ is good, if
$ \Gr^{F}M $ is coherent over
$ \Gr^{F}\cD_{X} $.
This condition implies that (1.1.1) holds for
$ r \gg 0 $ (i.e., the generating level is finite) locally on
$ X $.

(ii) The generating level of
$ (M,F) $ coincides with the minimum of the length of the filtrations of
filtered
$ \cO_{X} $-Modules
$ (L,F) $ having strict surjective morphisms
$ (\cD_{X},F)\otimes_{\cO_{X}}(L,F) \to (M,F) $,
where the filtration
$ F $ on the left-hand side is defined by
$ F_{p} = \sum_{i} F_{i}\cD_{X}\otimes_{\cO_{X}}F_{p-i}L $.
We say that
$ (\cD_{X},F)\otimes_{\cO_{X}}(L,F) $ is the induced filtered (left)
$ \cD_{X} $-Module associated with
$ (L,F) $.
The condition (1.1.1) is equivalent to the strict surjectivity of the
canonical morphism
$ (\cD_{X},F)\otimes_{\cO_{X}}(F_{r}M,F) \to (M,F) $.

(iii) Let
$ i : X \to Y $ be a closed embedding of complex manifolds, and
consider
$ i_{*}(M,F) $ as filtered
$ \cD $-Module.
Then
$ (M,F) $ has generating level
$ r $ if and only if
$ i_{*}(M,F) $ does, and
$ (M,F) $ is
$ r $-generated if and only if
$ i_{*}(M,F) $ is
$ (r-d) $-generated, where
$ d = \codim_{Y}X $.
See (1.5.3).

(iv) If
$ X = pt $,
and
$ (M,F) $ underlies a Hodge structure (where
$ F^{p} = F_{-p}) $,
the generating level of
$ (M,F) $ coincides with the level of Hodge structure which is defined by
$$
\max\{p \in \bZ : \Gr_{F}^{p}M \ne 0\} - \min\{p \in \bZ :
\Gr_{F}^{p}M \ne 0\}.
\leqno(1.1.2)
$$
If
$ \dim X \ge 1 $ and
$ (M,F) $ underlies a variation of Hodge structure with stalkwise level
$ r $,
then
$ (M,F) $ has generating level
$ \le r $,
but the equality does not hold in general.

(v) Let
$ 0 \to (M',F) \to (M,F) \to (M'',F)
\to $ be a short exact sequence of filtered $ \cD_{X}
$-Modules. If
$ (M',F), (M'',F) $ are
$ r $-generated, then
$ (M,F) $ is also
$ r $-generated.
If
$ (M,F) $ is
$ r $-generated, then
$ (M'',F) $ is
$ r $-generated, but
$ (M',F) $ is not necessarily.

\medskip\noindent
{\bf 1.2.}
Let
$ X $ be a complex manifold,
$ S $ a polydisc with a coordinate
$ t $,
and
$ X' = X\times S $ so that
$ X $ is identified with
$ X\times \{0\} $ in
$ X' $.
We denote by
$ V^{0}\cD_{X'} $ the subring of
$ \cD_{X'} $ generated by
$ \cO_{X'}, \cD_{X} $ and
$ t\partial_{t} $,
and define the decreasing filtration
$ V $ of
$ \cD_{X'} $ by
$$
V^{i}\cD_{X'} = t^{i}V^{0}\cD_{X'},\quad V^{-i}\cD_{X'} =
\sum_{0\le j\le i} \partial_{t}^{j}V^{0}\cD_{X'}
$$
for
$ i \ge 0 $.

We say that a coherent
$ \cD_{X'} $-Module
$ M $ admits the rational filtration
$ V $ along
$ X $ (see [16]),
if
$ M $ has an exhaustive decreasing filtration
$ V $ indexed by
$ \bQ $ such that

\smallskip
\noindent
(i)
$ V^{\alpha}M $ are coherent
$ V^{0}\cD_{X'} $-sub-Modules of
$ M $,

\noindent
(ii)
$ tV^{\alpha}M \subset V^{\alpha +1}M $ (with the equality for
$ \alpha > 0) $,

\noindent
(iii)
$ \partial_{t}V^{\alpha}M \subset V^{\alpha -1}M $,

\noindent
(iv)
$ \partial_{t}t - \alpha $ is nilpotent on
$ \Gr_{V}^{\alpha}M $ for
$ \alpha \in \bQ $.

\smallskip
\noindent
Here
$ \Gr_{V}^{\alpha}M = V^{\alpha}M/V^{>\alpha}M $ with
$ V^{>\alpha}M = \bigcup_{\beta >\alpha}V^{\beta}M $,
and we assume
$ V $ is indexed discretely (i.e., there exists a discrete subset
$ \Sigma $ of
$ \bR $ contained in
$ \bQ $ such that
$ V^{\alpha}M = V^{\beta}M $ if
$ \Sigma \cap [\alpha, \beta ) = \emptyset ) $.
The filtration
$ V $ is uniquely determined by the above conditions.
See also [8] [12].
If
$ M $ is a right
$ \cD_{X'} $-Module, we replace
$ \partial_{t}t - \alpha $ by
$ t\partial_{t} + \alpha $ in the above condition.
See (1.5.1).

We say that a coherent filtered
$ \cD_{X'} $-Module
$ (M,F) $ admits the rational filtration
$ V $ along
$ X $ (see [16]),
if
$ M $ admits the rational filtration
$ V $ along
$ X $ and the following conditions are satisfied:
$$
\alignedat2
t : (V^{\alpha}M, F) &\simto (V^{\alpha
+1}M, F)\quad
&&\text{for}\,\,\alpha > 0,
\\
\partial_{t} : \Gr_{V}^{\alpha}(M, F) &\simto
\Gr_{V}^{\alpha -1}(M, F[-1])\quad
&&\text{for}\,\,\alpha <1,
\endalignedat
\leqno(1.2.1)
$$
where
$ (F[m])_{p} = F_{p-m} $.

If the above conditions are satisfied, we define the nearby and vanishing
cycle functors
$ \psi_{t}, \varphi_{t} $ by
$$
\psi_{t,1}(M, F) = \Gr_{V}^{1}(M,F),\quad \varphi_{t,1}(M, F) =
\Gr_{V}^{0}(M,F[-1]),
$$
$$
\psi_{t,\ne 1}(M, F) = \varphi_{t,\ne 1}(M, F) = \bigoplus_{0<\alpha <1}
\Gr_{V}^{\alpha}(M,F),
$$
where
$ \psi_{t} = \psi_{t,1} \oplus \psi_{t,\ne 1} $ and the same for
$ \varphi_{t} $.
We have the morphisms
$$
\alignat2
\can : \psi_{t,1}(M, F)
&\to \varphi_{t,1}(M, F),\quad
&\Var : \varphi_{t,1}(M, F) &\to \psi_{t,1}(M, F[-1]),
\\
N : \psi_{t}(M, F)
&\to \psi_{t}(M, F[-1]),\quad
&N : \varphi_{t}(M, F) &\to \varphi_{t}(M, F[-1])
\endalignat
$$
by
$ - \partial_{t} $,
$ t $ and
$ - (\partial_{t}t - \alpha ) $ on
$ \Gr_{V}^{\alpha} $.
See [16].
For right
$ \cD $-Modules, we use (1.5) and shift the filtration
$ F $ by one so that
$ \varphi_{t} \circ (i_{0})_{*} = id $,
where
$ i_{0} : X \to X' $ denotes a natural inclusion.

\medskip\noindent
{\bf Remarks.} (i) The filtration
$ V $ is functorial and exact.
In fact, if
$ 0 \to M' \to M \to M'' \to 0 $ is a
short exact sequence of coherent
$ \cD_{X'} $-Modules such that
$ M $ admits the rational filtration
$ V $,
then
$ M', M'' $ admit the rational filtration
$ V $ and we have exact sequences
$ 0 \to V^{\alpha}M' \to V^{\alpha}M \to
V^{\alpha}M'' \to 0 $.
See [16, 3.1.5].

If
$ \supp M \subset X $,
we have
$ \Gr_{V}^{\alpha}M = 0 $ for
$ - \alpha \notin \bN $.
See [16, 3.1.3].
(We have also
$ \Gr_{V}^{0}{\circ}(i_{0})_{*} = id $ by definition of direct image.
See (1.5.3).)
So a morphism
$ M \to M' $ induces an isomorphism
$ V^{\alpha}M \to V^{\alpha}M' $ for
$ \alpha > 0 $ if
$ M|_{X'\setminus X} \to M'|_{X'\setminus X} $ is an
isomorphism.

(ii) If
$ (M,F) $ underlies a mixed Hodge Module
$ \cM $,
$ (M,F) $ admits the rational filtration
$ V $ and the conditions in (1.2.1) are satisfied.
Furthermore, the functors
$ \psi_{t} $,
$ \varphi_{t} $ and the morphisms
$ \can,
\Var, N $ are defined in the category of mixed Hodge Module in a
compatible way with the above definition, where the shift of
filtration is replaced with Tate twist. (For example, $ N $ is induced
by $ - (2\pi i)^{-1}\log T_{u} $,
where
$ T_{u} $ is the unipotent part of the monodromy
$ T $.)
See [16] [17].
In particular, the morphisms
$
\can,
\Var, N $ are strictly compatible with the Hodge filtration
$ F $.
See [16, 5.1.14].

(iii) If
$ \cM $ is a pure (i.e. polarizable) Hodge Module of weight
$ w $,
the weight filtration
$ W $ on
$ \psi_{t}\cM $ and
$ \varphi_{t,1}\cM $ is the monodromy filtration with center
$ r = w - 1 $ and
$ w $ respectively,
which is characterized by
$$
NW_{i} \subset W_{i-2},\quad N^{j} : \Gr_{r+j}^{W} \simto
\Gr_{r-j}^{W} (j > 0).
$$
Since
$ N $ is strictly compatible with
$ W $,
the functor assigning
$ \Gr_{i}^{W} $ commutes with
$ \Ker\, N, \Coker\, N $,
and
we have
$$
\Gr_{i}^{W}\Ker\, N = 0\,\,\,(i > r),\quad \Gr_{i}^{W}\Coker\, N =
0\,\,\,(i < r).
$$

(iv) Let
$ \cM $ be a polarizable Hodge Module on
$ X $,
and
$ Z $ an irreducible (or, more generally, pure dimensional) closed
analytic subset of
$ X $.
We say that
$ \cM $ has {\it strict support}
$ Z $,
if
$ \supp \cM = Z $ and
$ \cM $ has no nontrivial sub nor quotient object whose support has
dimension
$ < \dim Z $.
In this case, the same property holds for the underlying
$ \cD_{X} $-Module and perverse sheaf (in particular, the latter is an
intersection complex with local system coefficients).
Every polarizable Hodge Module
$ \cM $ has the strict support decomposition
$ \cM = \bigoplus_{Z} \cM_{Z} $ such that
$ \cM_{Z} $ has strict support
$ Z $.
See [16].

Let
$ \cM $ be a polarizable Hodge Module with strict support
$ Z $.
Then there exists a proper closed analytic subset
$ Z_{0} $ of
$ Z $ such that
$ Z \setminus Z_{0} $ is smooth and the restriction of
$ \cM $ to
$ Z \setminus Z_{0} $ is a polarizable variation of Hodge structure.
This means that the direct image of the variation of Hodge structure by
the closed embedding
$ Z \setminus Z_{0} \to X \setminus Z_{0} $ (whose underlying
filtered
$ \cD $-Module is defined in (1.5)) is isomorphic to
$ \cM|_{X\setminus Z_{0}} $.
Furthermore,
$ \cM $ is uniquely determined by the polarizable variation of Hodge
structure on
$ Z \setminus Z_{0} $,
and we have an equivalence of categories as in [17, 3.21].

Let
$ \cM $ be a polarizable Hodge Module with strict support
$ Z $.
We say that
$ \cM $ is of simple normal crossing type if
$ Z $ is smooth and the above
$ Z_{0} $ is a divisor
$ D $ with simple normal crossings.
Here simple means that the irreducible components of
$ D $ are smooth.

(v) Let
$ \cM $ be a polarizable Hodge Module of simple normal crossing
type with strict support $ X $,
and
$ D $ as above.
Let
$ (M,F) $ be the underlying filtered
$ \cD_{X} $-Module of
$ \cM $.
Then
$ (M,F)|_{X\setminus D} $ underlies a variation of Hodge structure.
Let
$ L^{\alpha} $ (resp.
$ L^{>\alpha}) $ be the Deligne extension of
$ M|_{X\setminus D} $ such that the eigenvalues of the residue of the
connection are contained in
$ [\alpha, \alpha + 1) $ (resp.
$ (\alpha, \alpha + 1]) $.
Then
$ L^{\alpha}, L^{>\alpha} $ are contained in
$ M(*D) $ the localization of
$ M $ along
$ D $.
By Schmid, the filtration
$ F $ can be extended to
$ L^{\alpha}, L^{>\alpha} $ so that
$ \Gr_{p}^{F}L^{\alpha}, \Gr_{p}^{F}L^{>\alpha} $ are locally free
$ \cO_{X} $-Modules.

Let
$ (x_{1}, \dots, x_{n}) $ be a local coordinate system of
$ X $ such that
$ D = \bigcup_{1\le i\le s} \{x_{i} = 0\} $ locally.
Let
$ V_{(i)} $ denote the rational filtration
$ V $ along
$ x_{i} = 0 $.
Then
$$
\aligned
L^{\alpha -1} &= L^{\alpha}\otimes_{\cO_{X}}\cO_{X}(D) \\
M \cap L^{\alpha -1} &= \bigcap_{1\le i\le s}
V_{(i)}^{\alpha}M \quad\text{for any}\,\, \alpha \in \bQ,\\
L^{\alpha -1} &= \bigcap_{1\le i\le s} V_{(i)}^{\alpha}M
\quad \text{for}\,\,
\alpha > 0
\endaligned
$$
(The last isomorphism means
$ L^{>-1} \subset M) $.
Using
$ L $,
the Hodge filtration
$ F $ on
$ M $ is expressed as:
$$
F_{p}M = \sum_{i} F_{i}\cD_{X}(F_{p-i}L^{>-1}).
\leqno(1.2.2)
$$
See [17, (3.10.12)].
In particular,
$ (M,F) $ has generating level
$ \le r $,
if the variation of Hodge structure on
$ X \setminus D $ has stalkwise level
$ r $.
We have also
$$
F_{p}M = \sum_{i} F_{i}\cD_{X}(M \cap F_{p-i}L^{-1}).
\leqno(1.2.3)
$$
This follows from the filtered isomorphisms for
$ \alpha < 1 $:
$$
\partial_{i} : \Gr_{V_{(i)}}^{\alpha}(M; F,V_{(j)} \,(j \ne i))
\simto \Gr_{V_{(i)}}^{\alpha -1}(M; F[-1],V_{(j)}\, (j \ne i))
\leqno(1.2.4)
$$
(see [17, 3.12]) if
$ M \cap F_{p-i}L^{-1} $ is replaced with
$ F_{p-i}M \cap L^{-1} $.
So the assertion is reduced to the strict injectivity of
$$
(M; F,V_{(j)}) \to (M(*D); F,V_{(j)}).
$$
See [17, 3.12] where
$ M, M(*D) $ are denoted by
$ {j}_{!*}^{\reg}M,\, {j}_{*}^{\reg}M $.

\medskip\noindent
{\bf 1.3.~Theorem.} {\it Let
$ f : X \to Y $ be a proper morphism of complex manifolds, and
$ Z $ an irreducible closed analytic subset of
$ X $.
Let
$ \cM $ be a polarizable Hodge Module on
$ X $ with strict support
$ Z $,
and
$ (M,F) $ its underlying filtered
$ \cD_{X} $-Module.
Assume
$ f $ is cohomologically
K\"ahler {\rm (}see {\rm [21]),} and the restriction of
$ \cM $ to
$ Z \setminus Z_{0} $ is a variation of Hodge structure whose stalkwise
level is
$ r $,
where
$ Z_{0} $ is a closed analytic subset of
$ Z $ containing
$ \Sing Z $.
Then
$ H^{j}f_{*}(M,F) $ {\rm (}see {\rm (1.5.5)} below{\rm )}
has generating level
$ \le r + \dim Z - |j| $ if
$ \dim f(Z) = 0 $,
and generating level
$ \le r + \dim Z - |j| - 1 $ otherwise.
}

\medskip\noindent
{\it Proof.}
Here we use right
$ \cD $-Modules.
See (1.5) below.
In particular,
$ (M,F) $ is a filtered right
$ \cD_{X} $-Module.
We first reduce the assertion to the case
$ X = Z $ and
$ \cM $ is a Hodge Module of simple normal crossing type.
See the above Remark (iv).
Let
$ \pi : X' \to Z $ be a projective morphism such that
$ X' $ is smooth,
$ D := \pi^{-1}(Z_{0}) $ is a divisor with simple normal crossings, and
$ \pi $ induces an isomorphism over
$ Z \setminus Z_{0} $,
where
$ Z_{0} $ is as above.
Let
$ \cM' $ be a polarizable Hodge Module with strict support
$ X' $ whose restriction to
$ X' \setminus D $ is isomorphic to that of
$ \cM $ to
$ Z \setminus Z_{0} $.
See [17, 3.21] (and the above Remark (iv)).
Let
$ (M',F) $ be the underlying filtered
$ \cD_{X'} $-Module of
$ \cM' $.
Then
$ (M,F) $ is a direct factor of
$ \pi_{*}(M',F) $ by the decomposition theorem for the underlying
filtered $ \cD $-Modules (see [21, (2.5)]).
So we may replace
$ X, \cM $ with
$ X', \cM' $,
and may assume
$ X = Z $ and
$ \cM $ is of simple normal crossing type
so that
$ D := Z_{0} $ is a divisor with simple normal crossings.
(Later
$ Z $ will be used to denote closed analytic subsets of
$ X $.)

Let
$ q = \min\{p \in \bZ : F_{p}M \ne 0\} $,
and
$ n = \dim X $.
Then
$$
\min\{p \in \bZ : F_{p}H^{j}f_{*}M \ne 0\} \ge q
$$
by definition of direct images for right
$ \cD $-Modules, and we have the hard Lefschetz property
$$
l^{ j} : H^{ -j}f_{*}(M,F) \simto H^{j}f_{*}(M,F[j])\quad
\text{for}\,\,j > 0,
$$
by shrinking
$ Y $ and replacing
$ X $ as above if necessary (where
$ (F[j])_{p} = F_{p-j}) $.
See [21].
So it is enough to show that
$ H^{j}f_{*}(M,F) $ is
$ (q + r + n) $- or
$ (q + r + n - 1) $-generated (depending on
$ \dim f(X)) $ because the hard Lefschetz property implies that
$ H^{j}f_{*}(M,F) $ is
$ (q + r + n - j) $- or
$ (q + r + n - j - 1) $-generated for
$ j > 0 $.

We consider first the case
$ \dim f(X) = 0 $.
Here we may assume
$ Y = pt $.
Then
$ f_{*}(M,F) $ is defined by
$ \bold{R}\Gamma (X, \DR_{X}(M,F)) $.
See (1.5) below.
Let
$ (x_{1}, \dots, x_{n}) $ be a local coordinate system such that
$ \bigcup_{i} \{x_{i} = 0\} \supset D $.
Let
$ V_{(i)} $ be as in the above Remark (v).
Then the
$ n+1 $ filtrations
$ F, V_{(i)} $ are compatible in the sense of [16].
Let
$ L $ denote the Deligne extension of the restriction to
$ X \setminus D $ of the left
$ \cD_{X} $-Module corresponding to
$ M $ such that the eigenvalues of the connection of the residue are
contained in
$ [0,1) $.
Then
$ {\Omega}_{X}^{n}\otimes_{\cO_{X}}L = \bigcap_{i} V_{(i)}^{1}M $.
See the above Remark (v).
Let
$ \DR_{X}(M)_{\log} $ be the intersection of
$ \DR_{X}(M) $ with the logarithmic
de Rham complex
$ \Omega_{X}(\log D)\otimes_{\cO_{X}}L[n] $.
Then (1.2.4) implies that the natural inclusion
$ \DR_{X}(M,F)_{\log} \to \DR_{X}(M,F) $ is a filtered
quasi-isomorphism, where the filtration
$ F $ on the left-hand side is the induced filtration.
Since
$ \Gr_{p}^{F}\DR_{X}(M)_{\log} = 0 $ for
$ p > q + r + n $ by (1.2.3), we get the assertion in this case.

In general, we proceed by induction on
$ \dim f(X) $.
Let
$ g $ be a function on an open subset
$ U $ of
$ Y $ such that
$ f^{ -1}g^{-1}(0) \ne f^{ -1}(U) $.
We may assume
$ f^{ -1}g^{-1}(0) \cup (D \cap f^{ -1}(U)) $ is a divisor with simple
normal crossings by the same argument as above.
By (1.4) below, it is enough to show that
$ \Gr_{V}^{\alpha}(i_{g})_{*}H^{j}f_{*}(M,F) $ is
$ (q + r + n - 1) $-generated for
$ 0 \le \alpha \le 1 $,
where
$ i_{g} $ is the embedding by the graph of
$ g $.
By [16, 3.3.17],
$ (f \times id)_{*}(i_{gf})_{*}(M; F,V) $ is strict for
$ (F, V) $ so that
$ \Gr_{V}^{\alpha} $ commutes with the cohomological direct image.
Since
$ i_{g^{\circ}}f = (f \times id){\circ}i_{gf} $,
we get
$$
\Gr_{V}^{\alpha}(i_{g})_{*}H^{j}f_{*}(M,F) =
H^{j}f_{*}\Gr_{V}^{\alpha}(i_{gf})_{*}(M,F).
$$
Let
$ (M^{\alpha},F) = \Gr_{V}^{\alpha}(i_{gf})_{*}(M,F) $.
We have the weight spectral sequence of filtered
$ \cD $-Modules
$$
{E}_{1}^{-k,k+j} = H^{j}f_{*}\Gr_{k}^{W}(M^{\alpha},F) \Rightarrow
H^{j}f_{*}(M^{\alpha},F)
$$
which underlies the weight spectral sequence of mixed Hodge Modules, and
degenerates at
$ E_{2} $.
Here the weight filtration
$ W $ is the monodromy filtration with center
$ w - 1 $ if
$ \alpha > 0 $,
and center
$ w $ if
$ \alpha = 0 $ (see Remark (iii) after (1.2)),
where
$ w $ is the weight of
$ \cM $.
So, using the semisimplicity of polarizable Hodge Modules, the assertion
is reduced to that
$ H^{j}f_{*}\Gr_{k}^{W}(M^{\alpha},F) $ is
$ (q + r + n - 1) $-generated for
$ 0 \le \alpha \le 1 $.
Here we may restrict to
$ 0 < \alpha \le 1 $,
because
$ t : \Gr_{k}^{W}(M^{0},F) \to \Gr_{k-2}^{W}(M^{1},F) $ is strictly
injective (see [16, 5.1.17]), and splits by the semisimplicity of
polarizable Hodge Modules.

We have the strict support decomposition
$$
\bigoplus_{0<\alpha \le 1} \Gr_{k}^{W}(M^{\alpha},F) = \bigoplus_{Z}
(M_{Z},F)
$$
where
$ Z $ are intersections of irreducible components of
$ D $,
and
$ (M_{Z},F) $ underlies a polarizable Hodge Module of simple normal
crossing type with strict support
$ Z $.
See the above Remark (iv).
So, using the inductive assumption, the assertion is reduced to
$ \Gr_{p}^{F}M_{Z}|_{Z\setminus Z_{0}} = 0 $ for
$ p > q + r + \codim_{X}Z - 1 $,
where
$ Z_{0} $ is the union of the intersection of
$ Z $ with the irreducible components of
$ D $ not containing
$ Z $.

Let
$ \Gr_{V}^{\nu} = \Gr_{V_{(1)}}^{\nu_{1}} \dots
\Gr_{V_{(n)}}^{\nu_{n}} $ for
$ \nu = (\nu_{1}, \dots, \nu_{n}) \in \bQ^{n} $,
where
$ V_{(i)} $ are as above.
Then the limit mixed Hodge structure of the variation of Hodge structure
$ (M_{Z},F)|_{Z\setminus Z_{0}} $ (up to the tensor with
$ {\Omega}_{Z}^{\dim Z} $ and the direct image by
$ Z \setminus Z_{0} \to X \setminus Z_{0}) $ is given by
$ \Gr_{V}^{\nu}(M_{Z},F) $ for
$ \nu \in (\bQ \cap [0,1])^{n} $ such that
$ \nu_{i} = 0 $ if and only if
$ Z \subset \{x_{i} = 0\} $.
So it is enough to show
$ \Gr_{p}^{F}\Gr_{V}^{\nu}M^{\alpha} = 0 $ for
$ \nu \in (\bQ \cap [0,1])^{n} $ and
$ p > q + r + s - 1 $,
where
$ s = \#\{i : \nu_{i} = 0\} $.
(Note that
$ \psi, \varphi $ in (1.2) induce exact functors of mixed Hodge Modules.)
By [17, (3.18.8)] we have
$$
\Gr_{V}^{\nu}(M^{\alpha},F) = \bigoplus_{0\le k<s'} \Gr_{V}^{\nu
+\alpha m}(M,F[k]),
$$
where
$ m = (m_{1}, \dots, m_{n}) $ is the multiplicity of
$ g $ (i.e.,
$ g = \prod_{i} x_{i}^{{m}_{i}} $ replacing
$ x_{i} $ if necessary), and
$ s' = \#\{i : \nu_{i} = 0, m_{i} \ne 0\} $.
So we get the assertion, because
$ \Gr_{p}^{F}\Gr_{V}^{\nu +\alpha m}M = 0 $ for
$ p > q + r $.

\medskip\noindent
{\bf Remarks.} (i) In the case
$ f $ is projective, we can use [16] instead of [21] in the proof of
(1.3).
Note that (0.1) is a special case of (1.3) where
$ f $ is the identity map.

(ii) Let
$ (M,F) $ be a coherent filtered (right)
$ \cD_{X} $-Module, and
$ (K^{\ssbull},F) \to (M,F) $ a filtered quasi-isomorphism of
filtered
$ \cD_{X} $-Modules such that each component
$ (K^{j},F) $ is the induced (right)
$ \cD_{X} $-Module associated with a filtered
$ \cO_{X} $-Module
$ (L^{j},F) $ (i.e.,
$ (K^{j},F) = (L^{j},F)\otimes_{\cO_{X}}(\cD_{X},F) $,
see [16, 2.1]).
Then
$ (L^{\ssbull},F) $ is a filtered differential complex in the sense of
[16, 2.2], and
$ \Gr_{p}^{F}L^{\ssbull} $ is a complex of
$ \cO_{X} $-Modules which is identified with
$ \Gr_{p}^{F}\Gr_{p}^{G}K $ where
$ G_{p}K^{j} = (F_{p}K^{j})\cD_{X} $ (see [loc. cit.]).
The level of
$ (K^{\ssbull},F) $ is defined by
$ \max \Lambda (K^{\ssbull},F) - \min \Lambda (K^{\ssbull},F) $,
where
$$
\Lambda (K^{\ssbull},F) = \{p \in \bZ : \Gr_{p}^{F}L^{j} \ne 0\,\,\,
\text{for some}\,\,j\}.
$$
Then the resolution level of
$ (M,F) $ is defined by the minimum of the level of
$ (K^{\ssbull},F) $ for filtered quasi-isomorphisms
$ (K^{\ssbull},F) \to (M,F) $ as above.
Clearly the resolution level of
$ (M,F) $ is greater than or equal to the generating level.
For example, the resolution level of
$ ({\omega}_{X},F) $ (with
$ \Gr_{p}^{F} = 0 $ for
$ p \ne 0) $ is
$ \dim X $,
and the generating level is
$ 0 $.

Similarly, we define the effective level of
$ (K^{\ssbull},F) $ and the effective resolution level of
$ (M,F) $ by using
$$
\Lambda '(K^{\ssbull},F) = \{p \in \bZ : H^{j}\Gr_{p}^{F}L^{\ssbull}
\ne 0\,\,\,
\text{for some}\,\,j\}.
$$
Note that the effective resolution level of
$ (K^{\ssbull},F) $ is independent of
$ (K^{\ssbull},F) $ which is filtered quasi-isomorphic to
$ (M,F) $ (using [16, 2.1.11] together with the equivalence of category
[16, 2.1.12]).
So the resolution level and the effective resolution level of
$ (M,F) $ coincide (taking a filtered quasi-isomorphism
$ (K^{\ssbull},F) \to (M,F) $ such that
$ K^{j} = 0 $ for
$ j > 0 $,
and using the above filtration
$ G) $.

Let
$ \DR_{X}(M,F) $ be the filtered
de Rham complex as in (1.5) below.
We define the
de Rham level of
$ (M,F) $ by
$ \max_{\Lambda''}(M,F) - \min_{\Lambda''}(M,F) $ with
$$
\Lambda ''(M,F) = \{p \in \bZ : H^{j}\Gr_{p}^{F}\DR_{X}(M)
\ne 0\,\,\, \text{for some}\,\,j\}.
$$
Then the
de Rham level coincides with the resolution level using the complex
of induced filtered $ \cD_{X} $-Modules associated with
$ \DR_{X}(M,F) $ (see [16, 2.2.6]).

Assume
$ (M,F) $ is Cohen-Macaulay (i.e.
$ \Gr_{F}M $ is Cohen-Macaulay over
$ \Gr_{F}\cD_{X}) $ so that the dual
$ \bD(M,F) $ (see [16, 2.4.3]) is isomorphic to a filtered
$ \cD_{X} $-Module
$ (M',F) $ up to a shift of complex.
Let
$$
q = \min\{p \in \bZ : F_{p}M \ne 0\},\quad q' = \min\{p \in \bZ :
F_{p}M' \ne 0\}.
$$
Then the resolution level of
$ (M,F) $ coincides with
$ - q' - q $ (using
$ \Gr_{-p}^{F}\bD(L^{\ssbull},F) = \bD\Gr_{p}^{F}(L^{\ssbull},F)
$,
where
$ \bD(L^{\ssbull},F) $ is the dual as filtered differential complex so
that
$ \bD(K^{\ssbull},F) $ is the complex of induced filtered
$ \cD_{X} $-Modules associated with
$ \bD(L^{\ssbull},F) $,
and
$ \bD $ on the right-hand side is the dual for
$ \cO_{X} $-Modules, see [loc. cit.]).

In particular, if
$ (M,F) $ underlies a polarizable Hodge Module with strict support
$ Z $ as in (1.3) such that the generic variation of Hodge structure has
stalkwise level
$ r $ and weight
$ w $,
then the weight of the Hodge Module is
$ w + m $ with
$ m = \dim Z $ so that
$ \bD(M,F) = (M,F[w + m]) $,
and the resolution level of
$ (M,F) $ is
$ r + m $,
because
$$
q = - (w + r)/2 - m,\quad q' = - (w + r)/2 + w.
$$
(Note that
$ F[n]_{p} = F_{p-n} $ and
$ F^{p} = F_{-p} $.)
As a corollary, we see that
$ (M,F) $ has generating level
$ \le r + m $.
This is slightly weaker than the assertion of (1.3).
(The latter implies that the generating level is strictly smaller than the
resolution level in this case.)

\medskip\noindent
{\bf 1.4.~Proposition.} {\it With the notation of (1.2), let
$ (M,F) $ be a coherent
$ \cD_{X'} $-Module admitting the rational filtration
$ V $ along
$ X, $ and assume
$ \Gr_{V}^{\alpha}(M,F) $ are
$ r $-generated {\rm (}as filtered
$ \cD_{X} $-Modules{\rm )} for
$ 0 \le \alpha \le 1 $.
Then the restriction of
$ (M,F) $ to an open neighborhood of
$ X $ is
$ r $-generated.
}

\medskip\noindent
{\it Proof.}
Since the filtration
$ V $ is indexed discretely, we see that
$ (V^{>0}M/V^{>1}M, F) $ is
$ r $-generated as filtered
$ \cD_{X} $-Modules.
Since
$ F_{p}(V^{>0}M/V^{>1}M) = F_{p}V^{>0}M/tF_{p}V^{>0}M $,
we get
$$
F_{i+r}V^{>0}M|_{X} = F_{i}V^{0}\cD_{X'}F_{r}V^{>0}M|_{X}\quad
\text{for}\,\,i > 0
\leqno(1.4.1)
$$
using Nakayama's lemma.
By the assumption for
$ \alpha = 0 $,
(1.4.1) holds with
$ V^{>0}M|_{X} $ replaced by
$ V^{0}M|_{X} $.
We have
$$
F_{i+r}V^{\alpha}M|_{X} \subset F_{i}\cD_{X'}F_{r}V^{>0}M|_{X}\quad
\text{for}\,\,i > 0
\leqno(1.4.2)
$$
by induction on
$ \alpha < 0 $ using the condition (1.2.1).
So we get the strict surjectivity of
$ (\cD_{X},F)\otimes_{\cO_{X}}(F_{r}M,F)|_{X} \to (M,F)|_{X} $.
(See Remark (ii) after (1.1).)
This implies the strict surjectivity of the morphism on an open
neighborhood of
$ X $ using the coherence of
$ (M,F) $.

\medskip\noindent
{\bf 1.5.}
In this paper we use mainly left
$ \cD $-Modules except in the proof of (1.3) and Remark (ii) after it.
Actually it is theoretically more natural to use right
$ \cD $-Modules in many places; for example, in the definition of
direct image below (see [16]), and the first equality of (0.4), (0.6)
(see [18]).

A filtered left
$ \cD_{X} $-Module
$ (M,F) $ corresponds to a right
$ \cD_{X} $-Module
$ (M^{r},F) $ so that
$$
(M^{r},F) = ({\Omega}_{X}^{n},F)\otimes_{\cO_{X}}(M,F),
$$
where
$ n = \dim X $.
If we choose a local coordinate system
$ (x_{1}, \dots, x_{n}) $,
then
$ {\Omega}_{X}^{n} $ is trivialized by
$ dx = dx_{1}\wedge \dots \wedge dx_{n} $ so that
$ M^{r} $ is identified with
$ M $ and the action of
$ \cD_{X} $ is given by using the involution
$ * $ of
$ \cD_{X} $ which is defined by
$$
(PQ)^{*} = Q^{*}P^{*},\,\,\,(x_{i})^{*} = x_{i},\,\,\,(\partial /\partial
x_{i})^{*}= - \partial /\partial x_{i}.
\leqno(1.5.1)
$$
Since the filtration
$ F $ on
$ {\Omega}_{X}^{n} $ is defined by
$ \Gr_{i}^{F}{\Omega}_{X}^{n} = 0 $ for
$ i \ne - n, $ we have
$$
F_{p}M^{r} = {\Omega}_{X}^{n}\otimes_{\cO_{X}}F_{p+n}M.
$$
The shift of filtration is necessary to get the isomorphism of the
de Rham complexes:
$ \DR_{X}(M,F) = \DR_{X}(M^{r},F) $,
where
$$
\aligned
F_{p}\DR_{X}(M)^{i} &= {\Omega}_{X}^{i+n}\otimes_{\cO_{X}}
F_{p+i+n}M, \\
F_{p}\DR_{X}(M^{r})^{i} &= \wedge^{-i}\Theta_{X}\otimes_{\cO_{X}}
F_{p+i}M^{r}.
\endaligned
\leqno(1.5.2)
$$
Let
$ f : X \to Y $ be a morphism of complex manifolds.
If
$ f $ is a closed embedding, we take a local coordinate system
$ (x_{1}, \dots, x_{m}) $ on
$ Y $ such that
$ X = \{x_{i} = 0 \,(i \le d)\} $.
Let
$ \partial_{i} = \partial /\partial x_{i} $.
Then we have locally
$$
f_{*}M = M\otimes_{\bC} \bC[\partial_{1}, \dots,
\partial_{d}],\quad F_{p}f_{*}M = \bigoplus_{\nu}
F_{p-|\nu |-d}M\otimes \partial^{\nu},
\leqno(1.5.3)
$$
where
$ \partial^{\nu} = \prod_{i} \partial_{i}^{\nu_{i}} $ for
$ \nu = (\nu_{1}, \dots, \nu_{d}) \in \bZ^{d} $.
This is compatible with the direct image for right
$ \cD $-Modules:
$$
f_{*}M^{r} = M^{r}\otimes_{\bC} \bC[\partial_{1}, \dots,
\partial_{d}],\quad F_{p}f_{*}M^{r} = \bigoplus_{\nu} F_{p-|\nu
|}M^{r}\otimes \partial^{\nu}.
\leqno(1.5.4)
$$
If
$ f $ is the projection
$ X = X_{0}\times Y \to Y $,
then
$ f_{*}(M,F) $ is a complex of filtered
$ \cD_{Y} $-Modules, and is defined by
$$
f_{*}(M,F) = \bold{R}f_{\ssbull}\DR_{X/Y}(M,F),
\leqno(1.5.5)
$$
where
$ \bold{R}f_{\ssbull} $ is the sheaf theoretic direct image, and
$ \DR_{X/Y} $ is defined by
$ \DR_{X_{0}} $ in (1.5.2).
We will denote by
$ H^{j}f_{*}(M,F) $ the cohomology of
$ f_{*}(M,F) $.
We have the same for right
$ \cD $-Modules.

\bigskip\bigskip\centerline{{\bf 2. Hypersurface Case}}

\bigskip\noindent
{\bf 2.1.}
Let
$ X, S, X' $ be as in (1.2).
Let
$ f : X \to S $ be a holomorphic function whose values are
contained in
$ S $,
and
$ i_{f} : X \to X' $ the embedding by graph of
$ f $.
Let
$$
(\cB_{f}, F) = (i_{f})_{*}(\cO_{X},F[-1])
$$
as filtered
$ \cD $-Module (see (1.5)), where the filtration
$ F $ of
$ \cO_{X} $ is defined by

\noindent

$ \Gr_{i}^{F}\cO_{X} = 0 $ for
$ i \ne 0 $.
Then we have
$$
\cB_{f} = \cO_{X}[\partial_{t}] \,
(= \cO_{X}\otimes_{\bC}\bC[\partial_{t}]),\quad
F_{p}\cB_{f} = \sum_{0\le i\le p} \cO_{X}\otimes
\partial_{t}^{i},
$$
so that the action of
$ \cD_{X'} $ is expressed by
$$
\xi (a \otimes \partial_{t}^{i}) = \xi a\otimes
\partial_{t}^{i} - (\xi f)a\otimes \partial_{t}^{i+1},\quad
t(a \otimes \partial_{t}^{i}) = fa\otimes \partial_{t}^{i}
- ia\otimes\partial_{t}^{i-1}
\leqno(2.1.1)
$$
for
$ \xi \in \Theta_{X}, a \in \cO_{X} $,
where the direct image
$ (i_{f})_{\ssbull} $ is omitted to simplify the notation.
See also [18] [19].
Let
$ \tilde\cB_{f} $ be the algebraic microlocalization of
$ \cB_{f} $ (see [19]) so that
$$
\tilde\cB_{f} = \cO_{X}[\partial_{t},\partial_{t}^{-1}],\quad
F_{p}\tilde\cB_{f} = \sum_{i\le p} \cO_{X}\otimes \partial_{t}^{i}.
$$
By [7] [8] [12],
$ \cB_{f} $ admits the rational filtration
$ V $ along
$ Y $,
and
$ \tilde\cB_{f} $ has the filtration
$ V $ in the sense of [19] such that
$ V^{\alpha}\tilde\cB_{f} $ are coherent
$ \cD_{X}[t, t\partial_{t}, \partial_{t}^{-1}] $-sub-Modules of
$ \tilde\cB_{f} $,

$ \partial_{t}t - \alpha $ is nilpotent on
$ \Gr_{V}^{\alpha}\tilde\cB_{f} $,
and
$$
t(V^{\alpha}\tilde\cB_{f}) \subset V^{\alpha
+1}\tilde\cB_{f},\quad \partial_{t} : V^{\alpha}\tilde\cB_{f}
\simto V^{\alpha -1}\tilde\cB_{f}.
\leqno(2.1.2)
$$
(Actually,
$ V $ is uniquely characterized by these conditions.)
By construction of
$ V $ in [19], we have
$$
V^{\alpha}\tilde\cB_{f} = \iota (V^{\alpha}\cB_{f}) +
F_{-1}\tilde\cB_{f}\quad
\text{for}\,\,\alpha \le 1,
\leqno(2.1.3)
$$
where
$ \iota : \cB_{f} \to \tilde\cB_{f} $ denotes a natural
inclusion.
In particular, we have
$$
\alignedat2
\iota^{-1}(V^{\alpha}\tilde\cB_{f}) &= V^{\alpha}\cB_{f}\quad
&&\text{for}\,\,\alpha \le 1,
\\
\iota : \Gr_{V}^{\alpha}\cB_{f} &\simto \Gr_{V}^{\alpha}\tilde\cB_{f}
\quad
&&\text{for}\,\,\alpha < 1.
\endalignedat
\leqno(2.1.4)
$$
As to the difference of
$ \iota^{-1}(V^{>1}\tilde\cB_{f}) $ and
$ V^{>1}\cB_{f} $,
we have
$$
\Ker(\Gr_{V}^{1}\cB_{f} \to \Gr_{V}^{1}\tilde\cB_{f})
= \Ker(N : \Gr_{V}^{1}\cB_{f} \to \Gr_{V}^{1}\cB_{f}),
\leqno(2.1.5)
$$
where
$ N $ is as in (1.2).
In fact,
$ N = \Var{\circ}
\can $
with the notation of (1.2), and
$ \can $
is identified with
$ \iota : \Gr_{V}^{1}\cB_{f} \to
\Gr_{V}^{1}\tilde\cB_{f} $ by (2.1.4).
Furthermore
$ \can $ is surjective and
$ \Var $ is injective.
See [16, 5.1.4].

\medskip\noindent
{\bf Remark.}
Let
$ b_{f}(s) $ be the
$ b $-function of
$ f $,
and
$ - \alpha_{f} $ the maximal root of
$ b_{f}(s)/(s+1) $.
See [18] [19].
Then
$ \alpha_{f} $ is positive by [7], and
$ F_{0}\cB_{f} \subset V^{>0}\cB_{f} $.
See [18, (1.7)].

\medskip\noindent
{\bf 2.2.}
For
$ X, f $ as above, let
$ M_{f}' = \cO_{X}[f^{ -1}], M_{f}'' = M_{f}'/\cO_{X} $
as in the introduction.
They underlie respectively the mixed Hodge Modules
$$
(j_{X\setminus Z})_{*}\bQ_{X\setminus Z}^{H}[n]
,\,\,\,(i_{Z})_{*}(i_{Z})^{!}\bQ_{X}^{H}[n+1],
$$
where
$ i_{Z} : Z \to X,\,
j_{X\setminus Z} : X \setminus Z
\to X $ are natural inclusions.
Furthermore the exact sequence
$ 0 \to \cO_{X} \to M_{f}' \to
M_{f}'' \to 0 $ underlies
$$
0 \to \bQ_{X}^{H}[n] \to (j_{X\setminus
Z})_{*}\bQ_{X\setminus Z}^{H}[n] \to
(i_{Z})_{*}(i_{Z})^{!}\bQ_{X}^{H}[n+1] \to 0.
\leqno(2.2.1)
$$
So
$ M_{f}', M_{f}'' $ have the Hodge filtrations
$ F $ such that
$ F $ on
$ M_{f}'' $ is the quotient filtration of
$ F $ on
$ M_{f}' $.

Let
$ M_{f} $ be as in the introduction so that
$ M_{f} $ underlies the pure Hodge Module
$ \IC_{Z}\bQ^{H} $ and has the Hodge filtration
$ F $.
See [16].
Then we have
$$
\Sol(M_{f}') = (j_{X\setminus Z})_{!}\bC_{X\setminus Z}[n],\quad
\Sol(M_{f}'') = \bC_{Z}[n-1],\quad \Sol(M_{f}) = \IC_{Z}\bC,
\leqno(2.2.2)
$$
where
$ \Sol(M_{f}') = \bold{R}{\cH}om_{\cD_{X}}(M_{f}',
\cO_{X}[n])\,(= \bD{\circ}\DR) $,
and
$ \IC_{Z}\bC $ is the intersection complex [1] (which is
the direct sum of
$ \IC_{Z_{i}}\bC $ where
$ Z_{i} $ are the irreducible components of
$ Z) $.
The natural morphism
$ \bC_{Z}[n-1] \to \IC_{Z}\bC $ is surjective in the
category of perverse sheaves (because
$ \IC_{Z}\bC $ has
non nontrivial quotient whose support has dimension
$ < n - 1) $,
and
$ M_{f} $ is a sub-Module of
$ M_{f}'' $ by the Riemann-Hilbert correspondence [9] [13].
See also [28].

In the level of mixed Hodge Modules,
$ \IC_{Z}\bQ^{H}(-1) $ is a subobject of
$ (i_{Z})_{*}(i_{Z})^{!}\bQ_{X}^{H}[n+1] $,
because
$$
\bD(\bQ_{X}^{H}[n]) = (\bQ_{X}^{H}[n])(n),\quad \bD(\IC_{Z}\bQ^{H}) =
\IC_{Z}\bQ^{H}(n-1).
$$
So
$ (M_{f},F[-1]) \to (M_{f}'',F) $ is strictly injective, and the
shift of filtration
$ F $ by one in (0.6), (0.8) comes from this.
(We do not have this shift of filtration if we use right
$ \cD $-Modules and replace
$ \cO_{X}, \bQ_{X}^{H}[n] $ with
$ \omega_{X}, \bD_{X}^{H}[-n] $,
where
$ \bD_{X}^{H} $ is the dual of
$ \bQ_{X}^{H}$.)

\medskip\noindent
{\bf Remark.}
With the notation of (2.1) and (2.2), we have an exact sequence of mixed
Hodge Modules
$$
0 \to \bQ_{Z}^{H}[n-1] \to \psi_{t,1}(i_{f})_{*}(\bQ_{X}^{H}[n])
\to \varphi_{t,1}(i_{f})_{*}(\bQ_{X}^{H}[n])\to 0,
$$
where
$ \bQ_{Z}^{H}[n-1] = (i_{Z})^{*}(\bQ_{X}^{H}[n]) $,
and the last morphism is
$ \can $ in (1.2).
See [17, 2.24].
So we get
$$
\bQ_{Z}^{H}[n-1] = \Ker\, N \subset \psi_{t,1}(i_{f})_{*}(\bQ_{X}^{H}[n])
\leqno(2.2.3)
$$
by the same argument as in (2.1.5).
We have furthermore
$$
\Gr_{i}^{W}(\bQ_{Z}^{H}[n-1]) = 0 \,(i \ge n)
,\,\,\,\Gr_{n-1}^{W}(\bQ_{Z}^{H}[n-1]) = \IC_{Z}\bQ^{H}.
\leqno(2.2.4)
$$
See [17, (4.5.7-9)].
Here
$ W $ is the weight filtration, and is compatible with (2.2.3).

\medskip\noindent
{\bf 2.3.}
With the notation of (2.1) and (2.2), we have a natural isomorphism
$$
\cB_{f}[t^{-1}] = (i_{f})_{*}M_{f}',
$$
and
$ \cB_{f}[t^{-1}] $ has the Hodge filtration
$ F $ such that
$$
(\cB_{f}[t^{-1}], F) = (i_{f})_{*}(M_{f}',F[-1]).
$$
This means that we have canonical isomorphisms
$$
\cB_{f}[t^{-1}] = M_{f}'[\partial_{t}],\quad F_{p}(\cB_{f}[t^{-1}]) =
\sum_{i\ge 0} F_{p-i}M_{f}'\otimes \partial_{t}^{i},
\leqno(2.3.1)
$$
where the action of
$ \cD_{X'} $ is given as in (2.1.1).
Since
$ \cB_{f}[t^{-1}] $ is the localization of
$ \cB_{f} $ by
$ t $,
$ \cB_{f}[t^{-1}] $ has also the rational filtration
$ V $ along
$ X $ such that
$ t : V^{\alpha}(\cB_{f}[t^{-1}]) \simto V^{\alpha +1}
(\cB_{f}[t^{-1}]) $.
Let
$ \iota ' : \cB_{f} \to \cB_{f}[t^{-1}] $ be a natural
morphism.
By Remark (i) after (1.2) we have
$$
\alignedat2
\iota^{\prime -1}(V^{\alpha}(\cB_{f}[t^{-1}])) = V^{\alpha}
\cB_{f}&\quad &&
\text{for}\,\,\alpha
\in \bQ,
\\
\iota ' : V^{\alpha}\cB_{f} \simto V^{\alpha}(\cB_{f}[t^{-1}])
&\quad &&\text{for}\,\,\alpha > 0.
\endalignedat
\leqno(2.3.2)
$$
By (2.3.1),
$ (\cB_{f}[t^{-1}], F) $ underlies a mixed Hodge Module so that the
conditions in (1.2.1) are satisfied and
$ t : \Gr_{V}^{0}(M_{f}, F) \simto \Gr_{V}^{1}(M_{f}, F) $.
See [17].
So we get
$$
\aligned
F_{p}(\cB_{f}[t^{-1}])& = \sum_{i\ge 0} \partial_{t}^{i}
F_{p-i}V^{0}(\cB_{f}[t^{-1}]),
\\
F_{p}V^{0}(\cB_{f}[t^{-1}]) &= V^{0}(\cB_{f}[t^{-1}]) \cap
j_{*}j^{*}F_{p}\cB_{f}
\endaligned
\leqno(2.3.3)
$$
by [16, 3.2.3],
where
$ j : X' \setminus X \to X' $ denotes a natural inclusion.

\medskip\noindent
{\bf Remark.}
By the same argument as in (2.1.5), we have
$$
\iota ' : \Gr_{V}^{0}\cB_{f} \simto \Im\,N \subset
\Gr_{V}^{0}(\cB_{f}[t^{-1}]),
\leqno(2.3.4)
$$
because
$ \Im\,N = \Im\,\can $
by the bijectivity of
$ \Var $ for
$ \cB_{f}[t^{-1}] $.

\medskip\noindent
{\bf 2.4.~Proof of (0.4) and (0.5).}
By (2.3.1) the natural projection
$$
\cB_{f}[t^{-1}] = M_{f}'[\partial_{t}] \to M_{f}'
\leqno(2.4.1)
$$
which send
$ \sum_{i} m_{i}\otimes \partial_{t}^{i} $ to
$ m_{0} $,
is strictly compatible with the Hodge filtration
$ F $.
So, by (2.3.3), it is enough to calculate the right-hand side of
$$
F_{p}V^{0}(\cB_{f}[t^{-1}]) = V^{0}(M_{f}'[\partial_{t}]) \cap
(\sum_{0\le i\le p} M_{f}'\otimes \partial_{t}^{i}).
\leqno(2.4.2)
$$
Since
$ t : V^{0}(M_{f}'[\partial_{t}]) \to V^{1}(M_{f}'[\partial_{t}]) $ and
$ t : \sum_{0\le i\le p} M_{f}'\otimes \partial_{t}^{i} \to
\sum_{0\le i\le p} M_{f}'\otimes \partial_{t}^{i} $ are bijective,
we may replace
$ V^{0} $ with
$ V^{1} $.
Here the inverse of the action of
$ t $ is given by
$$
t^{-1}(a \otimes \partial_{t}^{k}) = \sum_{0\le i\le k} (k!/i!)(a/f^{
k+1-i})\otimes \partial_{t}^{i}.
\leqno(2.4.3)
$$
Then we may replace
$ \cB_{f}[t^{-1}], M_{f}' $ with
$ \cB_{f}, \cO_{X} $.
See Remark (i) after (1.2).
By (2.1.2) and (2.1.4), we have
$ a\otimes 1 \in V^{k+1}\tilde\cB_{f} $ if and only if
$ a\otimes \partial_{t}^{k} \in V^{1}\cB_{f} $.
So we get the last assertion of (0.4) taking the composition of (2.4.3)
and (2.4.1), because they are
$ \cD_{X} $-linear.
The first assertion of (0.4) is clear by (2.3.3).

As for (0.5), the first condition is equivalent to the surjectivity of the
natural inclusion
$$
\omega_{X} \otimes_{\cO_{X}} F_{0}M_{f}'' \to \omega_{Z},
$$
(using Remark (i) below), where the morphism is defined by using the
Poincar\'e residue.
(In fact,
$ \omega_{Z} = \omega_{X}(Z)/\omega_{X} $,
where
$ \omega_{X}(Z) = f^{ -1}\omega_{X}) $.
So the condition is equivalent to
$ \tilde{V}^{1}\cO_{X} = \cO_{X} $ by the first equality of (0.4),
and hence to the second condition of (0.5) by (2.1.4) and [18, (1.7)].

\medskip\noindent
{\bf Remarks.} (i) Let
$ \bD(M_{f}'',F) $ denote the dual of
$ (M_{f}'',F) $ as complexes of filtered
$ \cD_{X} $-Modules.
See [16].
Then
$ \bD(M_{f}'',F) $ underlies
$ \bQ_{Z}^{H}(n)[n-1] $ (see (2.2)), and it is a filtered
$ \cD_{X} $-Module.
We can show that
$$
\DR_{X}(\bD(M_{f}'',F[n])) \,
(= \bD(\DR_{X}(M_{f}'',F[n])))
$$
is isomorphic to
$ (\underline{\Omega}_{Z},F) $ in the derived category of filtered
differential complexes.
Then, by the theory of duality (see also Remark (ii) after (1.3)), we get
$$
\omega_{X} \otimes_{\cO_{X}} F_{0}M_{f}'' = \Gr_{-n}^{F}
(\DR_{X}(M_{f}'')) = \bD(\Gr_{F}^{0}{\underline{\Omega}}_{Z}).
$$

(ii) We can also show that rational singularity is
du Bois in general.

\medskip\noindent
{\bf 2.5.~Proof of (0.6).}
The last assertion is reduced to the second and (0.4), because
$ (M_{f},F) \to (M_{f}'',F) $ is strictly injective.
For the second assertion, we have to show that
$ f^{ -k-1}\tilde{V}^{>k+1}\cO_{X} $ is annihilated by the composition
$ M_{f}' \to M_{f}'' \to M_{f}''/M_{f} $.
Let
$ m \in \tilde{V}^{>k+1}\cO_{X} $.
Using the direct image by
$ i_{f} $,
it is enough to show the vanishing of the image of
$ t^{-1}(m \otimes \partial_{t}^{k}) $ (see (2.4.3)) by the composition
$$
\cB_{f}[t^{-1}] \to \cB_{f}[t^{-1}]/\cB_{f} =
(i_{f})_{*}M_{f}'' \to (i_{f})_{*}(M_{f}''/M_{f})
$$
Since
$ \supp (i_{f})_{*}(M_{f}''/M_{f}) \subset X $,
we have
$ V^{>0}(i_{f})_{*}(M_{f}''/M_{f}) = 0 $.
See Remark (i) after (1.2).
By (2.4), we have
$ t^{-1}(m \otimes \partial_{t}^{k}) \in V^{0}(\cB_{f}[t^{-1}]) $.
So it is enough to show that the image in
$ \Gr_{V}^{0}(i_{f})_{*}(M_{f}''/M_{f}) $ is zero.

Since
$ m\otimes \partial_{t}^{k} \in \iota^{-1}(V^{>1}\tilde\cB_{f}) $ in the
notation of (2.1), its image in
$ \Gr_{V}^{1}\cB_{f} = \Gr_{V}^{1}(\cB_{f}[t^{-1}]) $ is contained in
$ \Ker\, N $.
See (2.1.5).
So the image of
$ t^{-1}(m \otimes \partial_{t}^{k}) $ in
$ \Gr_{V}^{0}(\cB_{f}[t^{-1}]) $ is also contained in
$ \Ker\, N $.
On the other side, the quotient
$ \Gr_{V}^{0}(i_{f})_{*}M_{f}'' $ of
$ \Gr_{V}^{0}(\cB_{f}[t^{-1}]) $ is identified with
$ \Coker\, N $.
See (2.3.4).
So it is enough to show the vanishing of the composition
$$
\Ker\, N \to \Coker\, N \to
\Gr_{V}^{0}(i_{f})_{*}(M_{f}''/M_{f}) = M_{f}''/M_{f}.
$$
See [16, 3.2.6] for the last isomorphism.
The first morphism is strictly compatible with
$ W $.
By Remark (iii) after (1.2) (with
$ r = n - 1) $,
we may replace
$ \Ker\, N $ with
$ \Gr_{n-1}^{W}\Ker\, N $ which is isomorphic to
$ M_{f} $ by (2.2.3-4).
So we get the second assertion, because
$ M_{f} $ has no nontrivial quotient whose support has dimension
$ < n - 1 $.

For the first assertion, assume
$ f^{ -1}(0) $ reduced as in the introduction.
Let
$ V $ denote also the filtration on
$ \cO_{X} $ induced by
$ V $ on
$ \cB_{f} $.
By (2.10) and (2.11) of [18], we have
$$
\gathered
V^{>1}\cO_{X} = f\cO_{X},\quad \omega_{Z} =
(\cO_{X}/f\cO_{X})\otimes_{\cO_{X}}\omega_{X},
\\
F_{1}M_{f}\otimes_{\cO_{X}}\omega_{X} = \pi_{*}\omega_{Z'}
= (\tilde{V}^{>1}\cO_{X}/f\cO_{X})\otimes_{\cO_{X}}
\omega_{X},
\endgathered
\leqno(2.5.1)
$$
where the tensor with
$ \omega_{X} $ comes from the transformation of left and right
$ \cD $-Modules in (1.5).
So we get the first assertion by multiplying the last term by
$ f^{ -1} $ which is induced by the inverse of the isomorphism
$ t : \Gr_{V}^{0}(\cB_{f}[t^{-1}]) \simto \Gr_{V}^{1}
(\cB_{f}[t^{-1}]) $ (see (2.4.3)).

\bigskip\bigskip\centerline{{\bf 3. Isolated Singularity Case}}

\bigskip\noindent
{\bf 3.1.}
With the notation of (2.1), let
$ Z = f^{ -1}(0) $,
and assume
$ n := \dim X > 1 $ and
$ \Sing Z = \{ 0 \} $
in this section.
We choose and fix a local coordinate system
$ (x_{1}, \dots, x_{n}) $ around
$ 0 $.

From now on,
$ \cO_{X,0} $ will be denoted by
$ A $ to simplify the notation.
Let
$ \partial_{i} = \partial /\partial x_{i}, f_{i} = \partial_{i}f $ so that
$ \{f_{i}\} $ is a regular sequence.
Then
$ \dim_{\bC}A/(\partial f) < \infty $,
where
$ (\partial f) $ is the ideal generated by
$ f_{i} \,(1 \le i \le n) $.
We define
$$
\overline{A}_{f} = A/\sum_{i\ne j} \Im(f_{i}\partial_{j} - f_{j}
\partial_{i}).
$$
Then
$ \overline{A}_{f} $ is isomorphic to Brieskorn's module
$ H_{f}'' := {\Omega}_{X,0}^{n}/df\wedge d{\Omega}_{X,0}^{n-2} $,
where
$ {\Omega}_{X,0}^{n} $ is trivialized by the local coordinates.
By [2] [22]
$ \overline{A}_{f} $ is a free
$ \bC\{t\} $-module of rank
$ \mu, $ where the action of
$ t $ is given by the multiplication by
$ f $.
It has also a meromorphic connection which is called the Gauss-Manin
connection.
In fact, the action of the inverse of
$ \partial_{t} $ on
$ \overline{A}_{f} $ is given by
$$
\partial_{t}^{-1}v = f_{i}u\quad
\text{with}\,\,\partial_{i}u = v,
\leqno(3.1.1)
$$
and is well-defined.

The localization
$ \overline{A}_{f}[t^{-1}] $ of
$ \overline{A}_{f} $ by
$ t $ is a regular holonomic
$ \cD $-module of one variable, and has the rational filtration
$ V $ which is characterized by the following conditions:
$ V^{\alpha} $ are finite over
$ \bC\{t\} $,
$ tV^{\alpha} \subset V^{\alpha +1} $ (with the equality for
$ \alpha \gg 0), \partial_{t}V^{\alpha} \subset V^{\alpha -1} $,
and
$ \partial_{t}t - \alpha $ is nilpotent on
$ \Gr_{V}^{\alpha}(\overline{A}_{f}[t^{-1}]) $.
See [16] [20].
We will denote also by
$ V $ the induced filtration on
$ \overline{A}_{f} $.
Then
$$
\alpha_{f} = \min\{\alpha \in \bQ : \Gr_{V}^{\alpha}\overline{A}_{f}
\ne 0\}
$$
by Remark (iii) below.
(See Remark after (2.1) for
$ \alpha_{f}$.)
It is called the minimal exponent (or the Arnold exponent) in this case.
We have
$$
\overline{A}_{f} \supset V^{>n-\alpha_{f}-1}(\overline{A}_{f}[t^{-1}]).
\leqno(3.1.2)
$$
See for example [20].

Let
$ \overline{A}_{f}[\partial_{t}] $ be the localization of
$ \overline{A}_{f} $ by
$ \partial_{t}^{-1} $.
It is also a regular holonomic
$ \cD $-module of one variable, and has the rational filtration
$ V $ as above.
Furthermore, the natural morphism
$ \overline{A}_{f}[\partial_{t}] \to \overline{A}_{f}[t^{-1}] $ is
strictly compatible with the filtration
$ V $,
and induces the isomorphisms
$$
V^{\alpha}(\overline{A}_{f}[\partial_{t}]) \simto
V^{\alpha}(\overline{A}_{f}[t^{-1}])\quad
\text{for}\,\,\alpha > 0.
$$
So they induce the same filtration
$ V $ on
$ \overline{A}_{f} $,
because
$ \overline{A}_{f} \subset V^{>0}(\overline{A}_{f}[\partial_{t}]) =
V^{>0}(\overline{A}_{f}[t^{-1}]) $.

\medskip\noindent
{\bf Remarks.} (i) We have natural isomorphisms
$$
\overline{A}_{f}[\partial_{t}] = H^{0}\DR_{X}(A[\partial_{t}]) =
H^{0}\DR_{X}(A[\partial_{t},\partial_{t}^{-1}]),
\leqno(3.1.3)
$$
where
$ A[\partial_{t}], A[\partial_{t},\partial_{t}^{-1}] $ are as in (2.1)
(with
$ A = \cO_{X,0}), $ and
$ \DR_{X} $ is the Koszul complex for
$ \partial_{1}, \dots, \partial_{n} $ shifted by
$ n $ in this case.
In fact, the first isomorphism follows from the theory of Gauss-Manin
system (see for example [15]), and the second from the bijectivity of the
action of
$ \partial_{t} $ on the middle term.
(A similar argument shows that
$ H^{j}\DR_{X}(A[\partial_{t}]) = H^{j}\DR_{X}(A[\partial_{t},
\partial_{t}^{-1}]) = 0 $ for
$ j \ne 1 - n, 0$.)
By (3.1.1) the isomorphisms are compatible with the action of
$ t, \partial_{t} $.

Furthermore, the isomorphisms in (3.1.3) are also compatible with the
filtration
$ V $.
Here the last two terms of (3.1.3) have the quotient filtration of
$ V $ on
$ A[\partial_{t}], A[\partial_{t},\partial_{t}^{-1}] $ in (2.1).
In fact, the assertion for the first isomorphism follows from [16, 3.4.8].
For the second, it is enough to verify the equality for
$ \alpha \ll 0 $,
and the assertion is reduced to
$ F_{0}(A[\partial_{t}]) \subset V^{>0}(A[\partial_{t}]) $ (see Remark
after (2.1)).

Let
$ \tilde{V} $ be the filtration on
$ A $ induced by the filtration
$ V $ on
$ A[\partial_{t},\partial_{t}^{-1}] $ as in the introduction.
Then its quotient filtration on
$ \overline{A}_{f} $ is contained in the filtration
$ V $ on
$ \overline{A}_{f} $ by the compatibility of (3.1.3) with
$ V $.

(ii) Let
$ H^{n-1}(X_{\infty}, \bC) $ denote the vanishing cohomology of
$ f $ at
$ 0 $ (i.e., the cohomology of the Milnor fiber), and
$ H^{n-1}(X_{\infty}, \bC)_{\lambda} = \Ker(T_{s} - \lambda ) $
for
$ \lambda \in \bC $,
where
$ T_{s} $ is the semisimple part of the monodromy
$ T $.
Let
$ e(\alpha) = \exp(2\pi i\alpha ) $.
Then we have isomorphisms
$$
\alignedat2
\Gr_{V}^{\alpha}(\overline{A}[t^{-1}]) &= H^{n-1}(X_{\infty},
\bC)_{e(-\alpha )}\quad &&
\text{for}\,\,\alpha \in \bQ,
\\
\Gr_{V}^{\alpha}\overline{A} &= H^{n-1}(X_{\infty},
\bC)_{e(- \alpha )}\quad &&
\text{for}\,\,\alpha > n - \alpha_{f} - 1,
\endalignedat
\leqno(3.1.4)
$$
where the last isomorphism follows from (3.1.2).

(iii) Let
$ {\tilde{H}}_{f}'' = \sum_{i\ge 0} (t\partial_{t})^{i}H_{f}''
$ the saturation of
$ H_{f}'' $.
Let
$ b_{f}(s) $ be the
$ b $-function of
$ f $ at
$ 0 $.
By [11]
$ b_{f}(s)/(s+1) $ is the minimal polynomial of the action of
$ - \partial_{t}t $ on
$ {\tilde{H}}_{f}'' / t{\tilde{H}}_{f}'' $.

\medskip\noindent
{\bf 3.2.}
Let
$ B = H_{[0]}^{n}\cO_{X} $ (the algebraic local cohomology).
See [9] [13].
Using the Cech cohomology,
$ B $ can be identified with
$ \bC[x_{1}^{-1}, \dots, x_{n}^{-1}]x^{-{\bold 1}} $ where
$ x^{-{\bold 1}} = (x_{1} \dots x_{n})^{-1} $,
and
$ \bold{1} = (1, \dots, 1) $.
So it is isomorphic to
$ \bC[\partial_{1}, \dots,\partial_{n}] $.
It is a unique simple regular holonomic
$ \cD_{X} $-Module supported on
$ \{0\} $.
We have a pairing of
$ A $ and
$ B $ by the composition
$$
A \times B \to B \to \bC,
\leqno(3.2.1)
$$
where the first morphism is by the action of
$ A $ on
$ B $,
and the second is the residue map which send
$ x^{-{\bold 1}} $ to
$ 1 $,
and
$ x^{-\nu -1} $ to
$ 0 $ for
$ \nu \ne 0 $.
In particular, the pairing is compatible with the action of
$ \cD_{X} $, i.e.,
$$
\langle P^{*}a, b \rangle = \langle a, Pb \rangle \quad
\text{for}\,\,P \in \cD_{X}, a \in A, b \in B,
$$
where
$ P^{*} $ is as in (1.5.1).
(Here it is more natural to put
$ A = \omega_{X,0} $,
and use the residue map
$ H_{[0]}^{n}\omega_{X} \to \bC $.)

\medskip\noindent
{\bf Remarks.} (i) Let
$ M_{f}, M_{f}', M_{f}'' $ be as in (2.2).
Let
$ j : Z \setminus \{0\} \to Z $ denote a natural inclusion.
Then we have
$$
\Sol(M_{f}) = \IC_{Z}\bC = \tau_{<0}\bold{R}j_{*}(\bC_{Z\setminus
\{0\}}[n-1]),
$$
where
$ \Sol $ is as in (2.2.2).
See [1]
for the last isomorphism.
So the exact sequence
$ 0 \to M_{f} \to M_{f}'' \to
M_{f}''/M_{f} \to 0 $ corresponds by the contravariant
functor
$ \Sol $ to the distinguished triangle
$$
\to R^{n-2}j_{*}\bC_{Z\setminus \{0\}} \to \bC_{Z}[n-1] \to
\tau_{<0}\bold{R}j_{*}(\bC_{Z\setminus \{0\}}[n-1]) \to.
$$
In particular,
$ R^{i}j_{*}\bC_{Z\setminus \{0\}} = H^{i}(Z\cap S_{\varepsilon},
\bC) = 0 $ for $ 0 < i < n-2 $ (as is well known, see [14]), where
$ S_{\varepsilon} $ is a sufficiently small sphere with center
$ 0 $ in
$ X $.

Let
$ H^{n-1}(X_{\infty}, \bC)^{T} $ be the invariant part of
$ H^{n-1}(X_{\infty}, \bC) $ by the action of the monodromy
$ T $.
Then we have natural isomorphisms
$$
H^{n-1}(X_{\infty}, \bC)^{T} = H^{n-2}(Z\cap S_{\varepsilon},
\bC) = \Hom_{\cD_{X}}(M_{f}', B).
\leqno(3.2.2)
$$
In fact the first isomorphism follows from
$$
H^{n-1}(X_{\infty}, \bC)^{T} = H^{n-1}(S_{\varepsilon}\setminus Z,
\bC) = H_{Z\cap S_{\varepsilon}}^{n}
(S_{\varepsilon}, \bC) =
H^{n-2}(Z\cap S_{\varepsilon}, \bC),
$$
where we use the Wang sequence for the first isomorphism.
For the second isomorphism of (3.2.2), we have
$$
\Hom_{\cD_{X}}(M_{f}''/M_{f}, B) =
\Hom_{\cD_{X}}(M_{f}'', B) =
\Hom_{\cD_{X}}(M_{f}', B),
$$
using
$ \Hom_{\cD_{X}}(M_{f}, B) = 0, \Hom_{\cD_{X}}(\cO_{X}, B) = 0 $.
So the assertion follows from the above distinguished triangle, because
$ M_{f}''/M_{f} $ is a direct sum of copies of
$ B $ and
$ \Sol(B) = \bC_{\{0\}} $.

(ii) By (2.2),
$ \IC_{Z}\bQ^{H}(-1) $ is a subobject of
$ (i_{Z})_{*}(i_{Z})^{!}\bQ_{X}^{H}[n+1] $,
and the support of the quotient is
$ \{0\} $ (shrinking
$ X $ if necessary).
So there exists a mixed
$ \bQ $-Hodge structure
$ H = (H_{\bC}, F, H_{\bQ}, W) $ with a short exact
sequence of mixed Hodge Modules on
$ X $:
$$
0 \to \IC_{Z}\bQ^{H}(-1) \to
(i_{Z})_{*}(i_{Z})^{!}\bQ_{X}^{H}[n+1] \to (i_{\{0\}})_{*}H
\to 0,
\leqno(3.2.3)
$$
where
$ H $ is identified with a mixed Hodge Module on
$ \{0\} $ (setting
$ F_{p} = F^{-p}) $,
and
$ i_{\{0\}} : \{0\} \to X $ is a natural inclusion.
In particular,
we have
$$
(M_{f}''/M_{f},F) = (i_{\{0\}})_{*}(H_{\bC},F)
\leqno(3.2.4)
$$
as filtered
$ \cD_{X} $-Modules, where
$ (i_{\{0\}})_{*} $ is as in (1.5).
Let
$$
r_{0} = \max\{p \in \bZ : \Gr_{F}^{p}H_{\bC} \ne 0\}.
$$
Then
$ (M_{f}''/M_{f},F) $ is exactly
$ (n-r_{0}) $-generated.
See (1.5.3).

By Remark (i) above,
$ H_{\bC} $ is the dual vector space of
$ H^{n-1}(X_{\infty}, \bC)^{T} $.
Using the theory of mixed Hodge Modules, we can express the mixed Hodge
structure
$ H $ as follows.
With the notation of (2.2), we have a distinguished triangle of mixed
Hodge Modules
$$
\to (j_{X\setminus Z})_{!}\bQ_{X\setminus Z}^{H}[n] \to
(j_{X\setminus Z})_{*}\bQ_{X\setminus Z}^{H}[n] \to
(i_{Z})_{*}(i_{Z})^{*}(j_{X\setminus Z})_{*}\bQ_{X\setminus Z}^{H}[n]
\to.
$$
See [6].
Applying the functor
$ \varphi_{f,1} $,
we have
$$
\to \psi_{f,1}\bQ_{X}^{H}[n]
\buildrel N \over \to \psi_{f,1}\bQ_{X}^{H}(-1)[n] \to
(i_{Z})^{*}(j_{X\setminus
Z})_{*}\bQ_{X\setminus Z}^{H}[n] \to.
$$
Then we take further
$ (i_{\{0\}})^{*} $,
and get
$$
\Coker(N : H^{n-1}(X_{\infty}, \bQ)_{1} \to
H^{n-1}(X_{\infty}, \bQ)_{1}(-1)) = H^{n}(i_{\{0\}})^{*}
(j_{X\setminus Z})_{*}\bQ_{X\setminus Z}^{H},
$$
using the isomorphism
$ H^{n-1}(X_{\infty}, \bQ)_{1} = (i_{\{0\}})^{*}\psi_{f,1}
\bQ_{X\setminus Z}^{H}[n] $.
(This is compatible with the mixed Hodge structure in [24].)

On the other hand, applying
$ (i_{0})^{*} $ to (3.2.3), we have
$$
H^{n+1}(i_{\{0\}})^{*}(i_{Z})^{!}\bQ_{X}^{H} = H,
$$
because
$ H^{j}(i_{\{0\}})^{*}\IC_{Z}\bQ = 0 $ for
$ j \ge 0 $.
The left hand-side is isomorphic to
$$
H^{n}(i_{\{0\}})^{*}(j_{X\setminus Z})_{*}\bQ_{X\setminus Z}^{H}
$$
using the distinguished triangle
$ \to (i_{Z})_{*}(i_{Z})^{!} \to id \to
(j_{X\setminus Z})_{*}(j_{X\setminus Z})^{*} \to $.
So we get
$$
H = \Coker(N : H^{n-1}(X_{\infty}, \bQ)_{1} \to H^{n-1}(X_{\infty},
\bQ)_{1}(-1))
\leqno(3.2.5)
$$
as mixed
$ \bQ $-Hodge structures.

\medskip\noindent
{\bf 3.3.~Proposition.} {\it With the notation of {\rm (3.1)} and
{\rm (3.2),} let
$ E_{f} = \Hom_{\cD_{X}}(M_{f}', B) $,
and
$ - \beta_{f} $ the minimal root of the
$ b $-function of
$ f $ at
$ 0 $.
Then, for
$ k > \beta_{f} - 1 $,
we have an injective morphism
$$
\phi_{k} : E_{f} \to B
\leqno(3.3.1)
$$
by assigning
$ u(f^{ -k}) $ to
$ u \in E_{f} $.
}

\medskip\noindent
{\it Proof.}
Let
$ b_{f}(s) $ be the
$ b $-function of
$ f $ so that
$$
b_{f}(s)f^{ s} = Pf^{ s+1}
\leqno(3.3.2)
$$
for
$ P \in \cD_{X,0}[s] $.
Then
$ M_{f}' $ is generated by
$ f^{ -k} $ for
$ k > \beta_{f} - 1 $,
substituting
$ s = - j - 1 $ to (3.3.2) for
$ j \ge k $.
So
$ \phi_{k} $ is injective for
$ k > \beta_{f} - 1 $.

\medskip\noindent
{\bf 3.4.~Remark.} By Remark (i) after (3.2) we have
$$
M_{f} = (\bigcap_{u\in E_{f}} \Ker\, u)/\cO_{X} \subset M_{f}''.
\leqno(3.4.1)
$$
This means for
$ a \in A $:
$$
a/f^k (\text{\rm mod}\,\, A) \in M_{f,0} \Leftrightarrow
a\phi_k(E_f) = 0.
\leqno(3.4.2)
$$
We may call this a trivial version of [28].
In general, it is not easy to determine
$ \phi_k(E_f) $,
except for the quasihomogeneous isolated singularity case
(see (4.4) below).

\bigskip\bigskip\centerline{{\bf 4. Quasihomogeneous Case}}

\bigskip\noindent
{\bf 4.1.}
With the notation and assumptions of (3.1), we assume further in this
section that
$ f $ is a quasihomogeneous polynomial of weight
$ w = (w_{1}, \dots, w_{n}) $,
i.e.,
$ f $ is a linear combination of monomials
$ x^{\nu} $ such that
$ \alpha_{w}(\nu) = 1 $,
where
$ w_{i} $ are rational numbers such that
$ 0 < w_{i} \le 1/2 $,
and
$ \alpha_{w}(\nu) = \sum_{i} w_{i}\nu_{i} $ for
$ \nu = (\nu_{1}, \dots, \nu_{n}) \in \bN^{n} $.

Let
$ \xi_{w} = \sum_{i} w_{i}x_{i}\partial_{i} $ so that
$ \xi_{w}f = f $.
With the notation of (3.1), (3.2), let
$$
A^{\alpha} = \Ker(\xi_{w}^{*} + \alpha ) \subset A,\quad B^{\alpha} =
\Ker(\xi_{w} + \alpha ) \subset B,
$$
where
$ \xi_{w}^{*} $ is as in (1.5.1).
(It is more natural to put
$ A = \omega_{X,0}$.)
Then
$ A^{\alpha}, B^{\alpha} $ are finite dimensional vector spaces on which
the pairing (3.2.1) induces a perfect pairing.
For example,
$ A^{\alpha} $ is spanned by the monomials
$ x^{\nu} $ such that
$ \alpha_{w}(\nu +{\bold 1}) = \alpha $,
where
$ {\bold 1} = (1, \dots, 1) $.
We have the (converging) infinite direct sum decomposition
$ A = \hat{\bigoplus}_{\alpha} A^{\alpha} $ and the direct sum
decomposition
$ B = \bigoplus_{\alpha} B^{\alpha} $.
We say that
$ a \in A $ is of degree
$ \alpha $ and
$ a\otimes \partial_{t}^{i} $ is of degree
$ \alpha - i $ if
$ a \in A^{\alpha} $.
We define the filtration
$ U $ on
$ A $ by
$$
U^{\alpha}A = \hat{\bigoplus}_{\beta \ge \alpha} A^{\beta}\,
(= A \cap\prod_{\beta \ge \alpha} A^{\beta})
$$
so that
$ \Gr_{U}^{\alpha}A = A^{\alpha} $.
The associated function
$ v_{w} $ on
$ A $ is defined by
$$
v_{w}(g) = \max\{\alpha \in \bQ : g \in U^{\alpha}A\}.
$$
These are naturally extended to
$ A[\partial_{t},\partial_{t}^{-1}] $ and
$ M_{f,0}' = A[f^{ -1}] $ by
$$
U^{\alpha}(A[\partial_{t},\partial_{t}^{-1}]) = \sum_{i}
(U^{\alpha +i}A\otimes \partial_{t}^{i}),\quad v_{w}(g \otimes
\partial_{t}^{i}) = v_{w}(g) - i,
$$
$$
U^{\alpha}M_{f,0}' = \bigcup_{k} f^{ -k}U^{\alpha +k}A,\quad v_{w}(g/f^{
k}) = v_{w}(g) - k.
$$
Since the morphisms
$ f_{i}\partial_{j} - f_{j}\partial_{i} $ are compatible with the above
direct sum decomposition, we get the direct sum decomposition
$ \overline{A}_{f} = \hat{\bigoplus}_{\alpha}
{\overline{A}}_{f}^{\alpha},$ such that
$$
{\overline{A}}_{f}^{\alpha} = A^{\alpha}/\sum_{i\ne j} \Im(f_{i}
\partial_{j} - f_{j}\partial_{i} : A^{\alpha +w_{i}+w_{j}-1} \to
A^{\alpha}).
$$
Then the filtration
$ U $ induces a filtration on
$ \overline{A}_{f} $ which is also denoted by
$ U $.
We define
$$
{\overline{B}}_{f}^{\alpha} = \bigcap_{i\ne j}
\Ker(f_{i}\partial_{j} - f_{j}\partial_{i} : B^{\alpha} \to
B^{\alpha +w_{i}+w_{j}-1}).
$$
Then the pairing (3.2.1) induces a perfect pairing of
$ {\overline{A}}_{f}^{\alpha} $ and
$ {\overline{B}}_{f}^{\alpha} $,
because
$ (f_{i}\partial_{j} - f_{j}\partial_{i} )^{*} =
- (f_{i}\partial_{j} - f_{j}\partial_{i}) $.

The ideal
$ (\partial f) $ (see (3.1)) is compatible with the direct sum
decomposition, and we have
$ A/(\partial f) = \bigoplus_{\alpha} (A/(\partial f))^{\alpha} $.

\medskip\noindent
{\bf Remarks.} (i) By Brieskorn (unpublished),
$ \overline{A}_{f} $ is stable by the action of
$ t\partial_{t} $,
and has a basis
$ \{v_{i}\} $ over
$ \bC\{t\} $ such that
$ \partial_{t}tv_{i} = \alpha_{i}v_{i} $ for
$ \alpha_{i} \in \bQ $ using a calculation like (4.1.4) below.
Here the
$ v_{i} $ are of degree
$ \alpha_{i} $,
and
gives a basis of
$ A/(\partial f) $.
So we have
$$
\sum_{1\le i\le \mu} \alpha_{i}t^{i} = (t^{w_{i}} - t)/(1 - t^{w_{i}})
\leqno(4.1.1)
$$
as is well known (using the morphism
$ (f_{1}, \dots, f_{n}) : \bC^{n} \to \bC^{n}) $.
Note that the left-hand side is the Poincar\'e polynomial of the graded
vector space
$ A/(\partial f) $.

The same calculation as above implies that the quotient filtration
$ U $ on
$ \overline{A}_{f} $ coincides with the filtration
$ V $ in (3.1).
They coincide further with the quotient filtration
$ \tilde{V} $ on
$ \overline{A}_{f} $ (see Remark (i) after (3.1)), because
$ U^{\alpha}A \subset \tilde{V}^{\alpha}A $ by the following Remark.

(ii) By an argument similar to [19, (3.2)], we have
$$
V^{\alpha}(A[\partial_{t},\partial_{t}^{-1}]) = \sum_{i}
\cD_{X,0}(U^{\alpha +i}A\otimes \partial_{t}^{i}).
\leqno(4.1.2)
$$
(In particular,
$ U^{\alpha}A \subset \tilde{V}^{\alpha}A $.)
In fact, let
$ {}'V^{\alpha} $ denote the right-hand side of (4.1.2).
We have the finiteness of
$ {}'V^{\alpha} $ over
$ \cD_{X,0}[\partial_{t}^{-1}] $ using the surjectivity of
$$
\sum f_{i} : \bigoplus U^{\alpha -1+w_{i}}A \to U^{\alpha}A\quad
\text{for}\,\,\alpha > n - \alpha_{f}.
\leqno(4.1.3)
$$
(For the last surjectivity we have
$ A^{\alpha} = (\partial f)^{\alpha} $ for
$ \alpha > n - \alpha_{f} $ by (4.1.1).)
Furthermore,
$ {}'V^{\alpha} $ is stable by the action of
$ t\partial_{t} $ because
$$
\xi_{w}(a \otimes 1) = \xi_{w}a\otimes 1 - \partial_{t}t(a \otimes 1)).
\leqno(4.1.4)
$$
This implies also that
$ \Gr_{{}'V}^{\alpha} $ is annihilated by
$ \partial_{t}t - \alpha $,
because we have for
$ a \in A^{\alpha} $
$$
(\xi_{w} + \alpha_{f})a = \alpha a,\quad (\xi_{w} + \alpha_{f})
(a \otimes 1) = -\xi_{w}^{*}(a \otimes 1) \in {{}'V}^{>\alpha}.
$$
So
$ {}'V $ satisfies the conditions of
$ V $ in (2.1).

(iii) By Remark (iii) after (3.1),
$ b_{f}(s)/(s+1) $ has only simple roots and the roots
are
$ \{- \alpha_{i}\} $.
In particular, we have by (4.1.1)
$$
\alpha_{f} = \alpha_{w}(1) = \sum_{i} w_{i},\quad \beta_{f} =
n - \alpha_{f}.
\leqno(4.1.5)
$$
See Remark after (2.1) for
$ \alpha_{f}, $ and (3.3) for
$ \beta_{f} $.

(iv) By the surjectivity of (4.1.3), we get the last equalities of (0.7),
(0.8),
because we have
$$
f^{ -k-1}(f_{i}U^{\alpha +w_{i}+k}A) \subset f^{ -k}U^{\alpha +k}A +
\partial_{i}(f^{ -k}U^{\alpha +w_{i}+k}A),
\leqno(4.1.6)
$$
using
$$
kf_{i}a/f^{ k+1} = (\partial_{i}a)/f^{ k} - \partial_{i}(a/f^{ k}).
\leqno(4.1.7)
$$

(v) Similarly
$ (M_{f}',F), (M_{f}'',F) $ have generating level
$ k_{0} $,
and
$ (M_{f},F) $ has generating level
$ k_{1} $ using (4.1.7) (together with the direct sum decomposition
$ A = \hat{\bigoplus}_{\alpha} A^{\alpha}) $ if we assume the first
equalities of (0.7), (0.8).
In fact, we have
$ A^{\beta_{f}} \ne (\partial f)^{\beta_{f}} $ by (4.1.1), (4.1.5), and
$ Af \subset (\partial f) $ by
$ \xi_{w}f = f $.

\medskip\noindent
{\bf 4.2.~Proof of (0.7).}
By Remarks (iv) and (v) after (4.1), it remains to show the first
equality.
By the same argument as in (2.4), it is enough to show
$$
F_{p}V^{\alpha}(A[\partial_{t},\partial_{t}^{-1}]) = \sum_{i\le p}
F_{p-i}\cD_{X,0}(U^{\alpha +i}A\otimes \partial_{t}^{i}).
\leqno(4.2.1)
$$
See (2.1) for
$ F $,
and (4.1.2) for
$ V $.
In fact, taking the intersection with
$ A[\partial_{t}] $ for
$ \alpha = 1 $,
we get
$$
F_{p}V^{1}(A[\partial_{t}]) = \sum_{0\le i\le p} F_{p-i}
\cD_{X,0}(U^{i+1}A\otimes \partial_{t}^{i}),
\leqno(4.2.2)
$$
using (2.1.4) for the left-hand side and
$ (\partial f) \subset U^{1}A $ (because
$ A = U^{\alpha_{f}}A) $ for the right.
Then the first equality follows by taking the composition of (2.4.3) and
(2.4.1).

By (4.1.2), (4.2.1) is reduced to
$$
(\sum_{i\le p} F_{p-i}\cD_{X,0}(U^{\alpha +i}A\otimes
\partial_{t}^{i})) \cap F_{p-1} = \sum_{i\le p-1} F_{p-1-i}
\cD_{X,0}(U^{\alpha +i}A\otimes \partial_{t}^{i}),
\leqno(4.2.3)
$$
using this formula inductively.
It is enough to show the inclusion
$ \subset $.
We define the filtration
$ V_{w} $ on
$ A, A[\partial_{t},\partial_{t}^{-1}] $ by
$$
V_{w}^{\alpha}A = U^{\alpha +\alpha_{f}}A,\quad
V_{w}^{\alpha}(A[\partial_{t},\partial_{t}^{-1}]) =
U^{\alpha +\alpha_{f}}(A[\partial_{t},\partial_{t}^{-1}]),
$$
(see also (5.3) below) so that
$ \Gr_{V_{w}}A := \bigoplus_{\alpha} \Gr_{V_{w}}^{\alpha}A $ is a graded
$ \bC $-algebra.
Then we may replace
$ U $ by
$ V_{w} $ in (4.2.3).
Take an element
$ m $ of the left-hand side.
Adding an element of the right-hand side if necessary, we may assume
$$
m = \sum_{i\le p} \sum_{|\nu |=p-i} \partial^{\nu}(a_{\nu}\otimes
\partial_{t}^{i})\quad
\text{with}\,\,a_{\nu} \in V_{w}^{\alpha +i}A,
$$
where
$ \partial^{\nu} = \prod_{i} \partial_{i}^{\nu_{i}} $ for
$ \nu = (\nu_{1}, \dots, \nu_{n}) $.
Let
$ {v}_{w}'(a) = v_{w}(a) - \alpha_{f} $ corresponding to
$ V_{w} $.
Let
$$
\alpha_{\nu} = {v}_{w}'(a_{\nu}),\quad \gamma_{\nu} =
\alpha_{\nu} + |\nu | - \alpha_{w}(\nu),
$$
so that
$ {v}_{w}'([\partial f]^{\nu}a_{\nu}) = \gamma_{\nu} $,
where
$ [\partial f]^{\nu} = \prod_{i} {f}_{i}^{\nu_{i}} $.
Let
$$
\aligned
\beta &= \min\{\gamma_{\nu} - p\},
\\
\Lambda (j)&= \{\nu \in \bN^{n} : \gamma_{\nu} - p = \beta,
|\nu | = p - j\},
\\
J &= \{ j \in J : \Lambda (j) \ne \emptyset \}.
\endaligned
$$
If
$ \beta \ge \alpha $,
$ m $ is contained in
$ \sum_{i\le p-1} V_{w}^{\alpha +i}A\otimes \partial_{t}^{i} $,
and the assertion is trivial.
So we may assume
$ \beta < \alpha $.
It is enough to show that, adding to
$ m $ an element of the right-hand side of (4.2.3),
$ m $ has an expression as above such that
$ \beta $ becomes larger.

We will first reduce to the case
$ \#J = 1 $.
Let
$ j' = \max J $.
Since
$ m $ belong to the left-hand side of (4.2.3), we have
$$
\sum_{j\in J}\sum_{\nu \in \Lambda (j)} (-1)^{|\nu |}[\Gr \,
\partial f]^{\nu}(\Gr \, a_{\nu}) = 0\quad
\text{in}\,\,\Gr_{V_{w}}^{\beta +p}A
\leqno(4.2.4)
$$
by (2.1.1), where
$ [\Gr \, \partial f]^{\nu} = \prod_{i} (\Gr \, f_{i})^{\nu_{i}} $.
This implies
$$
\sum_{\nu \in \Lambda (j')} [\Gr \, \partial f]^{\nu}(\Gr \, a_{\nu}) \in
(\Gr \, \partial f)^{p-j'+1},
$$
where
$ (\Gr \, \partial f) $ is the ideal of
$ \Gr_{V_{w}}A $ generated by
$ \Gr \, f_{i} \in \Gr_{V_{w}}^{1-{w}_{i}}A \,(1 \le i \le n) $.
Then we have
$ a'_{\nu,i} \in V_{w}^{{\alpha}_{\nu}+{w}_{i}-1}A $ for
$ \nu \in \Lambda (j') $ and
$ 1 \le i \le n $ such that
$$
\Gr \, a_{\nu} = \sum_{i} \Gr \, f_{i}a'_{\nu,i}\quad
\text{in}\,\,\Gr_{V_{w}}^{{\alpha}_{\nu}}A.
$$
See Remark (ii) below.
This implies
$$
\sum_{i} \Gr \, \partial_{i}(a'_{\nu,i}\otimes \partial_{t}^{j'-1}) =
(\sum_{i} (\Gr\,\partial_{i}a'_{\nu,i})\otimes \partial_{t}^{j'-1}
- \Gr \,a_{\nu}\otimes \partial_{t}^{j'}\quad
\text{in}\,\,\Gr_{V_{w}}^{{\alpha}_{\nu}j'}
([\partial_{t},\partial_{t}^{-1}]).
$$
So, replacing
$ m $ with
$ m - \sum_{\nu \in \Lambda(j')} \sum_{i} \partial^{\nu}
(\partial_{i}a'_{\nu,i}\otimes \partial_{t}^{j'-1}) $,
it has an expression as above such that
$ \beta $ does not decrease and
$ \max J - \min J $ becomes smaller.
Repeating the argument, the assertion is reduced to the case
$ J = \{j'\} $.

Let
$ A[T] $ be the polynomial ring over
$ A $ with variables
$ T = (T_{1}, \dots, T_{n}) $,
and define the filtration
$ V_{w} $ by
$$
V_{w}^{\alpha}(A[T]) = \sum_{\nu} V_{w}^{\alpha +{\alpha}_{w}
(\nu)-|\nu |} A T^{\nu},
$$
(i.e.,
$ T_{i} $ has degree
$ 1 - w_{i}) $ so that
$ \sum_{\nu} a_{\nu}\otimes T^{\nu} \in V_{w}^{\beta +p}(A[T]) $.

By (4.2.4) we have
$$
\sum_{\nu \in I(j')} \Gr \, a_{\nu}\otimes T^{\nu} = \sum_{i\ne k} \Gr \,
(f_{i}\otimes T_{k} - f_{k}\otimes T_{i})P_{i,k}(T)\quad
\text{in}\,\,\Gr_{V_{w}}^{\beta}A[T],
$$
where
$ P_{i,k}(T) \in V_{w}^{\beta -2+{w}_{i}+{w}_{k}}(A[T]) $ is homogeneous
in variables
$ T = (T_{1}, \dots, T_{n}) $ with degree
$ p - j' - 1 $.
See Remark (i) below.
So the assertion is reduced to the case where
$$
\sum_{\nu \in I(j')} \Gr \, a_{\nu}\otimes T^{\nu} = \Gr(f_{i}\otimes
T_{k} - f_{k}\otimes T_{i})(a \otimes T^{\nu})\quad
\text{in}\,\,\Gr_{V_{w}}^{\beta}A[T]
$$
for some
$ 0 < i < k \le n $ and
$ \nu \in \bN^{n} $ such that
$$
|\nu | = p - j' - 1,\quad {v}_{w}'(a) =
\beta + j' + \alpha_{w}(\nu) - 1 + w_{i} + w_{k}.
$$
In other words, we may assume
$$
m = \partial^{\nu}\partial_{i}(f_{k}a\otimes
\partial_{t}^{j'}) - \partial^{\nu}\partial_{k}(f_{i}a\otimes
\partial_{t}^{j'})
$$
for
$ i, k, \nu $ as above.
Since
$$
\partial_{i}(f_{k}a\otimes \partial_{t}^{j'}) - \partial_{k}
(f_{i}a\otimes \partial_{t}^{j'}) = \partial_{i}(\partial_{k}a
\otimes \partial_{t}^{j'-1}) - \partial_{k}(\partial_{i}a\otimes
\partial_{t}^{j'-1}),
$$
$ m $ belongs to the right-hand side of (4.3.2), and we get the assertion.

\medskip\noindent
{\bf Remarks.} (i) Let
$ R $ be a
$ \bC $-algebra with a morphism
$ \bC[y] \to R $,
where
$ \bC[y] $ is the polynomial ring over
$ \bC $ in variables
$ y = (y_{1}, \dots, y_{n}) $.
Let
$ g_{i} $ be the image of
$ y_{i} $ in
$ R $.
We assume
$ R $ is flat over
$ \bC[y] $.
(For example,
$ R = \Gr_{V_{w}}A $ with
$ g_{i} = \Gr \, f_{i}$.)
Let
$ R[T] $ be the polynomial ring over
$ R $ in variables
$ T = (T_{1}, \dots, T_{n}) $,
and
$ R[T]^{m} $ the
$ R $-submodule of
$ R[T] $ generated by homogeneous polynomials of degree
$ m $.
Then we have the exact sequence
$$
\bigoplus_{i\ne j} R[T]^{m-1} \to R[T]^{m} \to R,
$$
where the first morphism is defined by
$ \{P_{i,j}(T)\} \to \sum_{i\ne j} (g_{i}\otimes T_{j} - g_{j}\otimes
T_{i})P_{i,j}(T) $,
and the second by substituting
$ T_{i} = g_{i} $.
In fact, the assertion is reduced to the case
$ R = \bC[y] $ and
$ g_{i} = y_{i} $,
because the functor
$ \otimes_{\bC[y]} R $ is exact.
Then the proof if easy.
(This exactness seems to be known to some specialists.)

(ii) With the above notation and assumption, let
$ I $ be the ideal of
$ R $ generated by
$ g_{i} $,
and
$ \Lambda $ a subset of
$ \{\nu \in \bN^{n} : |\nu | = m\} $.
Then we have the exact sequence
$$
\bigoplus_{\nu \in \Lambda} I \to \bigoplus_{\nu \in \Lambda} R
\to R/I^{m+1},
$$
where the first morphism is a natural inclusion, and the second is induced
by the multiplication by
$ \prod_{i} {g}_{i}^{\nu_{i}} $.
This is also reduced to the case
$ R = \bC[y] $,
because
$ I = \Im(\sum_{i} g_{i} : \bigoplus_{i} R \to R) $.

\medskip\noindent
{\bf 4.3.~Proof of (0.8).}
The first inclusion follows from (0.6), because
$ U^{\alpha}A \subset \tilde{V}^{\alpha}A $ by (4.1.2).
(In this case, it is easy to show
$ u(a/f^{ k}) = 0 $ for
$ a \in U^{>k}A, u \in E_{f} $ using the action of
$ \xi_{w}$.)
So, by Remark (iv) and (v) after (4.1), it remains to prove the first
equality.
It is enough to show
$$
M_{f,0} \cap (\sum_{k\ge 0} F_{p-k}\cD_{X,0}G_{k}^{\ge
0}M_{f,0}'') \subset \sum_{k\ge 0}
F_{p-k}\cD_{X,0}G_{k}^{>0}M_{f,0}'',
\leqno(4.3.1)
$$
because the opposite inclusion is clear and the left-hand side coincides
with
$ F_{p+1}M_{f,0} $ using (0.7) and the strict injectivity of
$ (M_{f},F[-1]) \to (M_{f}'',F) $.
Take an element
$ m $ of the left-hand side.
It is represented by
$$
m' = \sum_{|\nu |\le p} \partial^{\nu}(a_{\nu}/f^{ |\nu |+1}) \in
M_{f}'\quad
\text{with}\,\,a_{\nu} \in U^{|\nu |+1}A,
$$
because
$ G_{k}^{\ge 0} \subset G_{k+1}^{\ge 0} $.
Since
$ m \in M_{f,0} $,
we get
$$
u(m') = \sum_{|\nu |\le p} \partial^{\nu}(a_{\nu}\phi_{|\nu
|+1}(u)) = 0\quad
\text{for any}\,\,u \in E_{f}.
$$
See (3.4.1).
We have
$ a_{\nu}\phi_{|\nu |+1}(u) \in B^{\alpha_{f}} $,
because
$ \phi_{|\nu |+1}(u) \in B^{|\nu |+1} $ (see (4.4) below) and
$ 1 \in A $ has degree
$ \alpha_{f} $.
Since
$ B $ is a free
$ \bC[\partial_{1}, \dots, \partial_{n}] $-module generated by
$ B^{\alpha_{f}} = \bC x^{-{\bold 1}} $ in the notation of (3.2),
we have $ a_{\nu}\phi_{|\nu |+1}(u) = 0 $,
and
$ \langle a_{\nu}, \phi_{|\nu |+1}(u)\rangle = 0 $ for any
$ \nu, u $.
Here the restriction of the paring (3.2.1) to
$ U^{\alpha}A \times U^{\alpha}B $ is factorized by
$ \Gr_{U}^{\alpha}A \times \Gr_{U}^{\alpha}B $.
So
$ \Gr_{U}^{|\nu |+1}a_{\nu} \in \Gr_{U}^{|\nu |+1}A $ belongs to
$ \sum_{i\ne j} \Im(f_{i}\partial_{j} - f_{j}\partial_{i}) $ by (4.4)
below, using the perfectness of the pairing between
$ {\overline{A}}_{f}^{|\nu |+1} $ and
$ {\overline{B}}_{f}^{|\nu |+1} $.
Then the assertion follows from
$$
\partial_{i}(\partial_{j}a/f^{ k}) - \partial_{j}(\partial_{i}a/f^{ k})
= - k(f_{i}(\partial_{j}a) - f_{j}(\partial_{i}a))/f^{ k+1},
$$
because
$ f_{i}(\partial_{j}a) - f_{j}(\partial_{i}a) \in U^{k+1}A $ and
$ \partial_{j}a \in U^{>k}A $ for
$ a \in U^{k + w_{i} + w_{j}}A $.

\medskip\noindent
{\bf Remarks.} (i) If
$ w_{i} = 1/d $,
the microlocal modified order (in the sense of [3]) of
$ a/f \in M_{f,0}' $ for
$ a \in A^{\alpha} $ is
$ d(\alpha - k) $ if
$ a/f \notin A $,
and
$ \infty $ otherwise (because the kernel of the natural morphism
$ M_{f,0}' \to \cE_{X,q}\otimes_{\cD_{X,q}}M_{f,0}' $ is
$ A $,
where
$ \cE_{X,q} $ is the ring of microdifferential operators at a general
point
$ q $ of
$ {T}_{0}^{*}X $,
see [9]).
So (0.8) implies a modified version (see [16]) of Brylinski's conjecture
using an argument as above.

(ii) By [23], we have
$$
\dim_{\bC}\Gr_{F}^{p}H^{n-1}(X_{\infty}, \bC)_{1} =
\dim_{\bC}(A/(\partial f))^{p}
$$
in the notation of Remark (ii) after (3.1) and (4.1).
(We can also use [27].)
By (3.2.5), the left-hand side coincides with
$ \dim_{\bC}\Gr_{F}^{p+1}H_{\bC} $ (because
$ N = 0) $.
So we get
$$
\dim_{\bC}\Gr_{F}^{p+1}H_{\bC} = \dim_{\bC}(A/(\partial
f))^{p}.
$$
If
$ f $ is homogeneous of degree
$ d $ (i.e.,
$ w_{i} = 1/d) $,
this formula implies that
$ r_{0} $ in Remark (ii) after (3.2) satisfies
$ r_{0} - 1 < n/d + 1 \le r_{0} $ using (4.1.1).
So
$ (M_{f}''/M_{f},F) $ is exactly
$ k_{0} $-generated, where
$ k_{0} $ is as in (0.7).
In particular, we see that
$ (M_{f}',F), (M_{f}'',F) $ are not
$ (k_{0}-1) $-generated.
See Remark (v) after (1.1).

If
$ f $ is not homogeneous, we can say only that
$ (M_{f}''/M_{f},F) $ is
$ k_{0} $-generated, because it is not easy to determine
$ r_{0} $ in general.

\medskip\noindent
{\bf 4.4.~Proposition.} {\it With the notation of {\rm (3.1)} and
{\rm (3.3),}
$ \phi_{k} $ induces the morphism
$$
\phi_{k} : E_{f} \to {\overline{B}}_{f}^{k},
\leqno(4.4.1)
$$
which is bijective for
$ k > \beta_{f} - 1 $,
and surjective for any
$ k $.
}

\medskip\noindent
{\it Proof.}
We have
$ u(f^{ -k}) \in {\overline{B}}_{f}^{k} $,
because
$ \xi_{w}(f^{ -k}) = - kf^{ -k} $ and
$ (f_{i}\partial_{j} - f_{j}\partial_{i})f^{ -k} = 0 $.
Since the monodromy
$ T $ is semisimple in the quasihomogeneous case,
$ H^{n-1}(X_{\infty}, \bC)_{1} $ coincides with the invariant part
$ H^{n-1}(X_{\infty}, \bC)^{T} $.
By the perfectness of the pairing, we have
$ \dim \overline{A}^{k} = \dim {\overline{B}}_{f}^{k} $,
and we get
$$
\dim {\overline{B}}_{f}^{k} = \dim E_{f}\quad
\text{for}\,\,k > n - \alpha_{f} - 1
\leqno(4.4.2)
$$
by (3.1.4), (3.2.2).
So we get the bijectivity by (3.3), because
$ \beta_{f} = n - \alpha_{f} $ (see (4.1.5)).
Then the surjectivity is reduced to the surjectivity of
$ f : {\overline{B}}_{f}^{k+1} \to {\overline{B}}_{f}^{k} $ which is
equivalent to the injectivity of
$ f : {\overline{A}}_{f}^{k} \to {\overline{A}}_{f}^{k+1} $ by the
perfectness of the pairing.
So the assertion follows from [22].

\medskip\noindent
{\bf Remark.}
By (3.4.1), (4.4) may be viewed as an algebraic version of [28].

\bigskip\bigskip\centerline{{\bf 5. Semiquasihomogeneous Case}}

\bigskip\noindent
{\bf 5.1.}
With the notation and the assumptions of (3.1), assume
$ f $ is semiquasihomogeneous of weight
$ w = (w_{1},\dots, w_{n}) $,
i.e.,
$ f = f' + f'' $,
where
$ f' $ is quasihomogeneous of weight
$ w $ with isolated singularity at
$ 0 $ (see (4.1)), and
$ f'' = \sum_{\nu} c_{\nu}x^{\nu} \in \bC\{x\} $ with
$ c_{\nu} = 0 $ for
$ \alpha_{w}(\nu) \le 1 $,
where
$ \alpha_{w} $ is as in (4.1).
Then we can define
$ \xi_{w}, A^{\alpha}, B^{\alpha}, U^{\alpha} $,
and
$ v_{w} $ as in (4.1) (but not
$ {\overline{A}}_{f}^{\alpha}, {\overline{B}}_{f}^{\alpha}) $.
We denote also by
$ U $ the quotient filtration on
$ \overline{A}_{f} $.
Then
$ \Gr_{U}^{\alpha}\overline{A}_{f} $ will be used later instead of
$ {\overline{A}}_{f}^{\alpha} $ in this section.

Let
$ K_{f} = \DR_{X}(A[\partial_{t},\partial_{t}^{-1}]) $,
the (algebraic) microlocal Gauss-Manin system.
See Remark (i) after (3.1).
Using the local coordinate system
$ (x_{1}, \dots, x_{n}) $,
$ K_{f} $ is identified with the Koszul complex (shifted by
$ n) $ for the action of
$ \partial_{1}, \dots, \partial_{n} $ on
$ A[\partial_{t},\partial_{t}^{-1}] $.
More precisely,
$ {K}_{f}^{j-n} $ is the direct sum of
$ A[\partial_{t},\partial_{t}^{-1}]\otimes dx_{p_{1}}\wedge \dots
\wedge dx_{p_{j}} $ for
$ 1 \le p_{1} < \dots < p_{j} \le n $.
The differential is defined by using (2.1.1).
We have
$ H^{j}K_{f} = 0 $ for
$ j \ne 1 - n, 0 $,
and
$ H^{0}K_{f} = \overline{A}_{f}[\partial_{t}] $.
See Remark (i) after (3.1).

We have the filtrations
$ F, U $ on
$ A[\partial_{t},\partial_{t}^{-1}] $ by
$ F_{p} = \bigoplus_{i\le p} A\otimes \partial_{t}^{i} $,
$ U^{\alpha} = \bigoplus_{i} U^{\alpha +i}A\otimes \partial_{t}^{i} $.
Then
$ K_{f} $ has the filtrations
$ F, U $ such that the restrictions of
$ F_{p}, U^{\alpha} $ to
$ A[\partial_{t},\partial_{t}^{-1}]\otimes dx_{p_{1}}\wedge \dots
\wedge dx_{p_{j}} $ are
$ F_{p+j-n}(A[\partial_{t},\partial_{t}^{-1}]) $ and
$ U^{\alpha +w'}(A[\partial_{t},\partial_{t}^{-1}]) $ respectively, where
$ w' = \alpha_{f'} - \sum_{k} w_{p_{k}} $ (i.e.
$ U $ is defined by counting also the weight of
$ dx_{i}) $.

\medskip\noindent
{\bf Remarks.} (i) The formula (4.1.2) remains true.
(In particular,
$ U^{\alpha}A \subset \tilde{V}^{\alpha}A $.)
Using
$ \Gr_{U} $,
the argument is almost the same as in Remark (ii) after (4.1).
In fact, the surjectivity of (4.1.3) follows from that of
$$
\sum_{i} \Gr \, f_{i} : \bigoplus_{i} \Gr_{U}^{\alpha -1+{w}_{i}}A \to
\Gr_{U}^{\alpha}A\quad
\text{for}\,\,\alpha > n - \alpha_{f'},
\leqno(5.1.1)
$$
combined with Nakayama's lemma.
(See (5.2.3) below for
$ \alpha_{f} = \alpha_{f'}$.)
Then we get the finiteness of
$ {}'V^{\alpha} $ over
$ \cD_{X,0}[\partial_{t}^{-1}] $.
Furthermore we have instead of (4.1.4)
$$
\xi_{w}(a \otimes 1) - \xi_{w}a\otimes 1 + \partial_{t}t(a \otimes 1))
\in {}'V^{>\alpha}
\leqno(5.1.2)
$$
for
$ a \in U^{\alpha}A $.
So the conditions of
$ V $ are verified.
See also [19, (3.2)].

(ii) As a corollary of the above remark, we get
$$
U = V\quad
\text{on}\,\,H^{0}K_{f} = \overline{A}_{f}[\partial_{t}],
\leqno(5.1.3)
$$
where
$ U $ on
$ H^{0}K_{f} $ is the quotient filtration of
$ U $ on
$ {K}_{f}^{0} = A[\partial_{t},\partial_{t}^{-1}] $,
$ V $ on
$ \overline{A}_{f}[\partial_{t}] $ is the filtration
$ V $ in (3.1),
and the last isomorphism of (5.1.3) is as in (3.1.1).
In fact,
$ V $ on
$ \overline{A}_{f}[\partial_{t}] $ is the quotient filtration of
$ V $ on
$ A[\partial_{t},\partial_{t}^{-1}] $ (see Remark (i) after (3.1)), and
$ H^{0}K_{f} = A[\partial_{t},\partial_{t}^{-1}]/\sum_{i} \Im\,
\partial_{i} $.
So (5.1.3) follows from (4.1.2) in this case.

\medskip\noindent
{\bf 5.2.~Proposition.} {\it The bifiltered complex
$ (K_{f}; F,U) $ in {\rm (5.1)} is strict in the sense of {\rm [16].}
}

\medskip\noindent
{\it Proof.}
By definition
$ \Gr_{p}^{F}K_{f} $ is the Koszul complex (shifted by
$ n) $ for the multiplication by
$ f_{1}, \dots, f_{n} $ on
$ A $,
and
$ H^{j}\Gr_{p}^{F}K_{f} = 0 $ for
$ j \ne 0 $.
This implies the injectivity of
$ H^{j}F_{p}K_{f} \to H^{j}F_{q}K_{f} $ for
$ p < q $ (i.e., the strictness of
$ F $ on
$ K_{f}) $ by using the long exact sequence associated with
$ 0 \to F_{p} \to F_{q} \to F_{q}/F_{p}
\to 0 $.

Similarly
$ \Gr_{p}^{F}\Gr_{U}K_{f} $ is the Koszul complex (shifted by
$ n) $ for the multiplication by
$ \Gr_{U}f_{1}, \dots, \Gr_{U}f_{n} $ on
$ \Gr_{U}A \,(= \bigoplus_{\alpha} \Gr_{U}^{\alpha}A) $,
and
$ H^{j}\Gr_{p}^{F}\Gr_{U}^{\alpha}K_{f} = 0 $ for
$ j \ne 0 $,
because
$ \Gr_{U}f_{i} = \Gr_{U}f'_{i} $.
This implies
$ H^{j}F_{p}\Gr_{U}^{\alpha}K_{f} = 0 $ for
$ j \ne 0 $,
using
$ F_{p}K_{f} \subset U^{\alpha}K_{f} $ for
$ \alpha + p \ll 0 $ (see Remark after (2.1)).
By the same argument as above, we get the injectivity of
$ H^{j}F_{p}U^{\alpha}K_{f} \to H^{j}F_{p}U^{\beta}K_{f} $ for
$ \alpha > \beta $.
So
$ H^{j}F_{p}U^{\alpha}K_{f} \to H^{j}K_{f} $ is injective (using
the strictness of
$ F) $.
In particular,
$ U $ on
$ K_{f} $ is strict.
By definition, the strictness of
$ (K_{f}; F,U) $ is equivalent to
$$
H^{j}F_{p}U^{\beta}K_{f} \cap H^{j}F_{q}U^{\alpha}K_{f} =
H^{j}F_{p}U^{\alpha}K_{f}\quad
\text{in}\,\,H^{j}K_{f}
$$
for
$ p < q, \alpha > \beta $.
This is verified by using the commutative diagram
$$
\CD
0 @>>> H^{j}F_{p}U^{\alpha}K_{f} @>>>
H^{j}F_{p}U^{\beta}K_{f} @>>>
H^{j}F_{p}(U^{\beta}/U^{\alpha})K_{f} @>>> 0 \\
@. @VVV @VVV @VVV @. \\
0 @>>> H^{j}F_{q}U^{\alpha}K_{f} @>>>
H^{j}F_{q}U^{\beta}K_{f} @>>>
H^{j}F_{q}(U^{\beta}/U^{\alpha})K_{f} @>>> 0.
\endCD
$$
where the injectivity of the right vertical morphism follows from
the vanishing of
$ H^{j}(F_{q}/F_{p})(U^{\beta}/U^{\alpha})K_{f} $ for
$ j \ne 0 $.
See [16].

\medskip\noindent
{\bf Remarks.} (i) As a corollary of (5.2), we have
$$
U = V\quad
\text{on}\,\,\overline{A}_{f},
\leqno(5.2.1)
$$
where
$ U $ is the quotient filtration of
$ U $ on
$ A $,
and
$ V $ is the induced filtration of
$ V $ on
$ \overline{A}_{f}[\partial_{t}] $.
In fact, we have
$ \overline{A}_{f} = F_{0}H^{0}K_{f} $ (because
$ \overline{A}_{f} $ is the image of
$ A $ in
$ H^{0}K_{f} $,
and is stable by
$ \partial_{t}^{-1}) $,
and
$ V $ on
$ \overline{A}_{f} $ is identified with
$ U $ on
$ F_{0}H^{0}K_{f} $ by (5.1.3).
By (5.2) we have the strict surjectivity of
$ ({K}_{f}^{0}; F,U) \to (K_{f}H^{0}; F,U) $ which implies that of
$ (F_{0}{K}_{f}^{0},U) \to (H^{0}F_{0}K_{f},U) $.
So we get
$$
V^{\alpha}\overline{A}_{f} = \sum_{k\ge 0} \partial_{t}^{-i}
U^{\alpha - i}\overline{A}_{f},
$$
and the assertion is reduced to
$ \partial_{t}^{-1}U^{\alpha -1}\overline{A}_{f} \subset
U^{\alpha}\overline{A}_{f} $.
But this is easily verified using (3.1.1).

(ii) We have a natural isomorphism
$$
H^{0}K_{f} = H^{0}K_{f'},\quad \Gr_{U}^{\alpha}\overline{A}_{f} =
\Gr_{U}^{\alpha}\overline{A}_{f'}.
\leqno(5.2.2)
$$
In fact, we have
$ \Gr_{U}(K_{f},F) = \Gr_{U}(K_{f'},F), $ and
$ \Gr_{U}H^{0}K_{f} = \Gr_{U}H^{0}K_{f'} $.
This implies the first isomorphism, because the filtration
$ U = V $ splits by the action of
$ \partial_{t}t $.
For the second isomorphism, we use
$$
\Gr_{U}^{\alpha}\overline{A}_{f} = \Gr_{U}^{\alpha}F_{0}H^{0}K_{f} =
H^{0}F_{0}\Gr_{U}^{\alpha}K_{f},
$$
where the last isomorphism follows from (5.2).
See [16].
(Note that the first isomorphism of (5.2.2) is related with [10], and is
not compatible with
$ F $.)

By [20], (5.2.2) implies that the exponents for
$ f $ and
$ f' $ coincide, and we get
$$
\alpha_{f} = \alpha_{f'} = \alpha_{w}(1) = \sum_{i} w_{i},\quad
\beta_{f} \le n - \alpha_{f}.
\leqno(5.2.3)
$$

(iii) The last isomorphisms of (0.7), (0.8) are true by the same argument
as in Remark (iv) after (4.1), because (4.1.3), (4.1.6) remain valid.

(iv) As a corollary of (5.2.2), we have the strictness of
$$
\sum_{i\ne j} (f_{i}\partial_{j} - f_{j}\partial_{i}) : \bigoplus_{i\ne j}
(A/U^{>\alpha +w_{i}+w_{j}-1}A, U[w_{i}+w_{j}-1]) \to
(A/U^{>\alpha}A, U)
\leqno(5.2.4)
$$
for
$ \alpha \in \bQ $,
where
$ (U[\beta ])^{\alpha} = U^{\alpha +\beta} $.
In fact, let
$ L^{-1} = \bigoplus_{i\ne j} (A, U[w_{i}+w_{j}-1]), L^{0} = A $,
and
$ d_{f} = \sum_{i\ne j} f_{i}\partial_{j} - f_{j}\partial_{i} : L^{-1}
\to L^{0} $ so that
$ \overline{A}_{f} = \Coker\, d_{f} $.
By (5.2.2) we have
$$
\dim_{\bC} \Im(\Gr_{U}d_{f} : \Gr_{U}^{\alpha}L^{-1} \to
\Gr_{U}^{\alpha}L^{0}) = \dim_{\bC} \Gr_{U}^{\alpha}
\Im\,d_{f},
$$
because this holds for
$ f' $ and
$ \Gr_{U}d_{f} = \Gr_{U}d_{f'} $.
Then the strictness of (5.2.4) follows.

\medskip\noindent
{\bf 5.3.}
We define the filtration
$ V_{w} $ on
$ \cD_{X,0} $ so that
$ V_{w}^{\alpha}\cD_{X,0} $ is the
$ \cO_{X,0} $-submodule of
$ \cD_{X,0} $ generated by
$ x^{\nu}\partial^{\nu '} $ for
$ \alpha_{w}(\nu) - \alpha_{w}(\nu ') \ge \alpha $.
Let
$ V_{w} $ be the filtration on
$ A $ so that
$ V_{w}^{\alpha}A $ is the ideal generated by
$ x^{\nu} $ for
$ \alpha_{w}(\nu) \ge \alpha $ (i.e.,
$ V_{w}^{\alpha}A = U^{\alpha +\alpha_{f}}A $,
see (4.1)).
We define the filtration
$ V_{w} $ on
$ M_{f,0}' = \cO_{X,0}[f^{ -1}] $ by
$ V_{w}^{\alpha}M_{f,0}' = \sum_{k\ge 0} f^{ -k}V_{w}^{\alpha
+k}A $ (i.e.,
$ V_{w}^{\alpha}M_{f,0}' = U^{\alpha +\alpha_{f}}M_{f,0}' $ if
$ f = f' $,
see (4.1)).
Then
$ (M_{f,0}',V_{w}) $ is a filtered
$ (\cD_{X,0}, V_{w}) $-Module, i.e.,
$ V_{w}^{\alpha}\cD_{X,0}V_{w}^{\beta}M_{f,0}' \subset
V_{w}^{\alpha +\beta}M_{f,0}' $.

Let
$ c $ be the smallest positive integer such that
$ {w}_{i}^{-1}c \in \bN $.
Then
$ V_{w} $ satisfies the following conditions:

\noindent
(i) $ V_{w}^{\alpha}M_{f,0}' $ are finitely generated
$ V_{w}^{0}\cD_{X,0} $-modules,

\noindent
(ii) $ V_{w}^{c}\cD_{X,0}V_{w}^{\alpha}M_{f,0}' =
V_{w}^{c+\alpha}M_{f,0}' $ for
$ \alpha \gg 0 $,

\noindent
(iii) $ \xi_{w} - \alpha $ is nilpotent on
$ \Gr_{V_{w}}^{\alpha}M_{f,0}' $.

(See also [8] in the case
$ w_{i} = 1/d $.)

In fact, we have the surjectivity of
$ U^{c}A \times U^{\alpha}A \to U^{\alpha +c}A $ for
$ \alpha \gg 0 $ (because
$ \alpha_{w}(\nu) \gg 0 $ implies
$ w_{i}\nu_{i} \ge k $ for some
$ i) $.
So (ii) follows.
Let
$ r $ be the smallest positive integer such that
$ c^{-1}rw_{i} \in \bN $ for any
$ i $.
Since
$ \{\Gr_{V_{w}}(f_{i})^{r}\} $ is a regular sequence, we have the
surjectivity of
$$
\sum_{i} (f_{i})^{r} : \bigoplus_{i} U^{\alpha -r+rw_{i}}A \to
U^{\alpha}A\quad
\text{for}\,\,\alpha \gg 0.
$$
Then
$$
f^{ -k}V_{w}^{\alpha +k}A \subset \sum_{i}
V_{w}^{r{w}_{i}}A\partial_{i}^{r}(f^{ -k+r}V_{w}^{\alpha +k-r}A)
+ \sum_{0\le j<k} f^{ -j}V_{w}^{\alpha +j}A\quad
\text{for}\,\,k \gg 0,
$$
and we get (i).
For (iii) we apply an argument similar to Remark (i) after (5.1).

We also define the filtration
$ V_{w} $ on
$ B $ so that
$ V_{w}^{\alpha}B $ consists of linear combinations of
$ x^{-\nu -1} (\nu \in \bN^{n}) $ with
$ \alpha_{w}(- \nu - 1) \ge \alpha $ in the notation of (3.2) and (4.1).

\medskip\noindent
{\bf Remarks.} (i) The filtration
$ V_{w} $ on
$ M_{f,0}' $ is uniquely characterized by the above three conditions.
In fact, if
$ V', V'' $ are two filtrations satisfying the conditions, take
$ \alpha \gg 0 $ such that the condition (ii) is satisfied for
$ V' $.
Then
$ V^{\prime \alpha} \subset V^{\prime \prime \beta} $ for some
$ \beta \le \alpha $,
and
$ V^{\prime \alpha +ic} \subset V^{\prime \prime \alpha} $ for
$ i \in \bN $ such that
$ \beta + ic \ge \alpha $.
So the image of
$ V^{\prime \alpha}/V^{\prime \alpha +ic} $ in
$ V^{\prime \prime \beta}/V^{\prime \prime \alpha} $ is annihilated by
$ P(\xi_{w}) $ where
$ P $ is a polynomial in one variable whose roots are greater than or
equal to
$ \alpha $.
This implies that the image is zero, and we get
$ V^{\prime \alpha} \subset V^{\prime \prime \alpha} $ for
$ \alpha \gg 0 $.
Then we can apply a similar argument for any
$ \alpha $,
and get the inclusion for any
$ \alpha \in \bQ $.
So the uniqueness follows.
(See also [8] [16].)

(ii) For
$ B $,
the filtration
$ V_{w} $ is characterized by the two conditions (i) and (iii).
In fact,
$ V_{w}^{\alpha}B $ are artinian
$ A $-modules by (i), and
$ \bigcap_{\alpha} V_{w}^{\alpha}B = 0 $,
because
$ B $ is simple.
Then
$ V_{w}^{\alpha}B = 0 $ for
$ \alpha \gg 0, $ and the condition (ii) is not necessary.

(iii) As a corollary of the above remarks, any
$ \cD_{X,0} $-linear morphism
$ M_{f,0}' \to B $ is strictly compatible with
$ V_{w} $,
because the quotient filtration on the image satisfies the conditions (i)
and (iii).

(iv) Using the local coordinate system
$ (x_{1}, \dots, x_{n}) $,
we have the inclusion
$ \Gr \, A := \Gr_{V_{w}}A \,(= \bigoplus_{\alpha} \Gr_{V_{w}}^{\alpha}A)
\to A $.
This induces the functor
$$
M \to \Sp_{w}(M) := A\otimes_{\Gr \, A}\Gr_{V_{w}}M,
\leqno(5.3.1)
$$
where
$ M $ is a finite
$ \cD_{X,0} $-module having the filtration
$ V_{w} $ satisfying the above three conditions.
We have natural isomorphisms as
$ \cD_{X,0} $-modules
$ \Sp_{w}(M_{f,0}') = M_{f',0}', \Sp_{w}(B) = B $.
So we get a map
$$
E_{f} \to E_{f'}.
\leqno(5.3.2)
$$
This is bijective, because it is injective by Remark (iii) above and
$ \dim E_{f} = \dim E_{f'} $ by (3.2.2) (combined with (3.1.4), (5.1.3),
(5.2.1)).

\medskip\noindent
{\bf 5.4.~Proof of (0.9).}
Except for the assertions on generating levels, the argument is
essentially the same as in the quasihomogeneous case.
In fact, the arguments in (4.2) and Remark (iv) after (4.1) hold also in
the semiquasihomogeneous case, and we get the (partial) generalization
of (0.7).
As to (0.8), we use the functor
$ \Sp_{w} $ in (5.3.1) and surjectivity of (5.3.2) to prove
$$
\Gr_{U}^{|\nu |+1}a_{\nu} \in \sum_{i\ne j}
\Im\,\Gr(f_{i}\partial_{j} - f_{j}\partial_{i})
$$
as in (4.3).
Then the assertion is reduced to the strictness of (5.2.4) for
$ \alpha = k + 1 $.

\medskip\noindent
{\bf Remarks.} (i) The assertions on generating levels in (0.7), (0.8) are
not true in the semiquasihomogeneous case.
For example, consider
$ f = x_{1}^{6} + x_{2}^{4} + x_{3}^{4} + x_{4}^{4} +
x_{1}^{2}x_{2}x_{3}x_{4} $.
Here
$ n = 4, \alpha_{f} = 11/12 $,
and
$ k_{0} = k_{1} = 2 $,
but
$ (M_{f},F), (M_{f}',F) $,
and
$ (M_{f}'',F) $ have generating level
$ \le 1 $ using (4.1.6).
In fact, we have
$$
U^{3}A = U^{3}(\partial f) + \bC
x_{1}^{4}x_{2}^{2}x_{3}^{2}x_{4}^{2},
$$
(same for
$ U^{>3}) $,
and
$$
-\xi_{w}^{*}(x_{1}^{2}x_{2}x_{3}x_{4}/f^{2}) =
- {1 \over 6}x_{1}^{4}x_{2}^{2}x_{3}^{2}x_{4}^{2}/f^{3},
$$
where
$ x_{i}x_{1}^{2}x_{2}x_{3}x_{4} \in U^{>2}A $.
In this case we have
$$
\Gr_{F}^{p}H^{3}(X_{\infty}, \bC)_{1}= 0\quad
\text{if and only if\quad} p \ne 2,
$$
and
$ r_{0} $ in Remark (ii) after (3.2) is
$ 3 $.
See Remark (ii) after (4.3).

(ii) If
$ w_{i} = 1/d $,
the assertion on generating level in (0.7) is generalized to the
semiquasihomogeneous case by Remark (i) after (4.3), because
$ \dim_{\bC}\Gr_{F}^{p}H^{n-1}(X_{\infty}, \bC)_{1} $ is
constant under a
$ \mu $-constant deformation.
See for example [26].
However, the assertion on generating level in (0.8) cannot be generalized
even in this case (e.g.,
$ f = x_{1}^{4} + x_{2}^{4} + x_{3}^{4} + x_{1}^{2}x_{2}^{2}x_{3}) $.

\bigskip\bigskip
\centerline{{\bf References}}

\bigskip

\item{[1]} Beilinson, A., Bernstein, J. and Deligne, P., Faisceaux
pervers, Ast\'erisque, vol. 100, Soc. Math. France, Paris, 1982.

\item{[2]} Brieskorn, E., Die Monodromie der isolierten
Singularit\"aten von Hyperfl\"achen, Manu\-scripta Math., 2 (1970),
103--161.

\item{[3]} Brylinski, J.-L., Modules holonomes \`a singularit\'es
r\'eguli\`eres et filtrations de Hodge II, Ast\'erisque vol. 101--102
(1983), 75--117.

\item{[4]} Deligne, P., Equation diff\'erentielle \`a points
singuliers r\'eguliers, Lect. Notes in Math. vol. 163, Springer,
Berlin, 1970.

\item{[5]} du Bois, Ph, Complex de de Rham filtr\'e d'une vari\'et\'e
singuli\`ere, Bull. Soc. Math. France, 109 (1981), 41--81.

\item{[6]} Durfee, A. and Saito, M., Mixed Hodge structure on the
intersection cohomology of links, Compos. Math., 76 (1990), 49--67.

\item{[7]} Kashiwara, M., B-function and holonomic systems, Inv. Math.
38 (1976), 33--53.

\item{[8]} \SameAuthor, Vanishing cycle sheaves and holonomic systems
of
differential equations, Lect. Notes in Math., vol. 1016, Springer,
Berlin, $ 1983, pp. 136-142 $.

\item{[9]} Kashiwara, M. and Kawai, T., On the holonomic system of
microdifferential equations III, Publ. RIMS, Kyoto Univ. 17 (1981),
813--979.

\item{[10]} L\^e, D. T. and Ramanujam, C. P., The invariance of Milnor
number implies the invariance of the topological type, Am. J. Math. 98
(1976), 67--78.

\item{[11]} Malgrange, B., Le polyn\^ome de Bernstein d'une
singularit\'e isol\'ee, in Lect. Notes in Math., vol. 459, Springer,
Berlin, 1975, pp. 98--119.

\item{[12]} \SameAuthor, Polyn\^ome de Bernstein-Sato et cohomologie
\'evanescente, Ast\'erisque, 101--102 (1983), 243--267.

\item{[13]} Mebkhout, Z., Une autre \'equivalence de cat\'egories,
Compo. Math. 51 (1984), 63--88.

\item{[14]} Milnor, J., Singular points of complex hypersurfaces, Ann.
Math. Stud. vol. 61, Princeton Univ. Press, 1969.

\item{[15]} Pham, F., Singularit\'es des systemes diff\'erentiels de
Gauss-Manin, Prog. in Math. vol. 2, Birkh\"auser, Boston, 1979.

\item{[16]} Saito, M., Modules de Hodge polarisables, Publ. RIMS,
Kyoto Univ., 24 (1988), 849--995.

\item{[17]} \SameAuthor, Mixed Hodge Modules, Publ. RIMS Kyoto Univ.
26 (1990), 221--333.

\item{[18]} \SameAuthor, On b-function, spectrum and rational
singularity, Math. Ann. 295 (1993), 51--74.

\item{[19]} \SameAuthor, On microlocal b-function, Bull. Soc. Math.
France 122 (1994), 163--184.

\item{[20]} \SameAuthor, On the structure of Brieskorn lattice, Ann.
Institut Fourier 39 (1989), 27--72.

\item{[21]} \SameAuthor, Decomposition theorem for proper K\"ahler
morphisms, T\^ohoku Math. J. 42 (1990), 127--148.

\item{[22]} Sebastiani, M., Preuve d'une conjecture de Brieskorn,
Manuscripta Math. 2 (1970), 301--308.

\item{[23]} Steenbrink, J, Intersection form for quasi-homogeneous
singularities, Compos. Math. 34 (1977), 211--223.

\item{[24]} \SameAuthor, Mixed Hodge structure on the vanishing
cohomology, in Real and Complex Singularities (Proc. Nordic Summer
School, Oslo, 1976) Alphen a/d Rijn: Sijthoff \& Noordhoff 1977, pp.
525--563.

\item{[25]} \SameAuthor, Mixed Hodge structures associated with
isolated singularities, Proc. Symp. Pure Math. 40 (1983) Part 2,
513--536.

\item{[26]} Varchenko, A., The complex exponent of a singularity does
not change along strata $ \mu $ = const, Func. Anal. Appl. 16 (1982),
1--9.

\item{[27]} \SameAuthor, The asymptotics of holomorphic forms
determine a mixed Hodge structure, Soviet Math. Dokl. 22 (1980)
772--775.

\item{[28]} Vilonen, K., Intersection homology $ D $-module on local
complete intersections with isolated singularities, Inv. Math. 81
(1985), 107--114.

\bye